\documentclass{llncs}

\usepackage{a4wide}
\usepackage{amsfonts}
\usepackage{latexsym}
\usepackage{graphicx}
\usepackage{epsfig}
\usepackage{amsbsy}
\usepackage{amsmath}
\usepackage{amssymb}
 \usepackage{graphicx,color}
 \usepackage{epsfig}
 \usepackage{curves,eepic}
 \usepackage{bm}
 \usepackage{bbm}
 \usepackage{soul}
 \usepackage{appendix}
 \usepackage{graphicx,rotating,booktabs}
\usepackage[margin=0.5in]{geometry}
\usepackage{caption}
\usepackage{multicol}

\newcommand{\beq}{\begin{equation}}
\newcommand{\deq}{\end{equation}}
\newcommand{\baq}{\begin{eqnarray}}
\newcommand{\daq}{\end{eqnarray}}

\newcommand{\baqm}{\begin{eqnarray*}}
\newcommand{\daqm}{\end{eqnarray*}}

\def\E{\mathbb{E}}

\usepackage[T1]{fontenc}
\usepackage[utf8]{inputenc}
\usepackage{tabularx,ragged2e,booktabs}
\newcolumntype{C}[1]{>{\Centering}m{#1}}

  \makeatletter
\def\Ddots{\mathinner{\mkern1mu\raise\p@
\vbox{\kern7\p@\hbox{.}}\mkern2mu
\raise4\p@\hbox{.}\mkern2mu\raise7\p@\hbox{.}\mkern1mu}}
\makeatother

\title{Heavy traffic analysis of a polling model with retrials and glue periods}
\author{Murtuza Ali Abidini$^*$, Jan-Pieter Dorsman$^{**}$ \and Jacques Resing$^*$}
\institute{$^*$
EURANDOM and Department of Mathematics and Computer Science\\
Eindhoven University of Technology\\
P.O. Box 513, 5600 MB Eindhoven, The Netherlands \\
\vspace{0.3cm}
$^{**}$ Korteweg–de Vries Institute for Mathematics\\
University of Amsterdam \\
P.O.Box 94248, 1090 GE Amsterdam, The Netherlands \\
\vspace{0.3cm}
\email{m.a.abidini@tue.nl,j.l.dorsman@uva.nl,j.a.c.resing@tue.nl}}
\begin{document}

\maketitle

\begin{abstract}
We present a heavy traffic analysis of a single-server polling model,
with the special features of retrials and glue periods. The combination of these features in a polling model typically occurs in certain optical networking models, and in models where
customers have a reservation period just before their service period. Just before the server arrives at a station there is some deterministic glue period.
Customers (both new arrivals and retrials) arriving at the station
during this glue period will be served during the visit of the server.
Customers arriving in
any other period leave immediately and will retry after an exponentially
distributed time. As this model defies a closed-form expression for the queue length distributions, our main focus is on their heavy-traffic asymptotics, both at embedded time points (beginnings of glue periods,
visit periods and switch periods) and at arbitrary time points. We obtain closed-form expressions for the limiting scaled joint queue length distribution in heavy traffic and use these to accurately approximate the mean number of customers in the system under different loads.
\end{abstract}
\numberwithin{equation}{section}

\section{Introduction}

Polling models are queueing models in which a single server, alternatingly,
visits a number of queues in some prescribed order. These models
have been extensively studied in the literature. For example, various
different service disciplines (rules which describe the server's behaviour
while visiting a queue) and both models with and without switch-over times
have been considered. We refer to Takagi \cite{Takagi1,Takagi2} and Vishnevskii and Semenova \cite{Vishnevskii}
for some literature reviews and to Boon, van der Mei and Winands \cite{Boon},
Levy and Sidi \cite{Levy} and Takagi \cite{Takagi3} for overviews of the
applicability of polling models.

Motivated by questions regarding the performance modelling and analysis
of optical networks, the study of polling models with {\it retrials} and
{\it glue periods} was initiated in Boxma and Resing \cite{BoxmaResing}.
In these models, just before the server arrives at a station there is
some glue period. Customers (both new arrivals and retrials) arriving at
the station during this glue period "stick" and will be served during the
visit of the server. Customers arriving in any other period leave
immediately and will retry after an exponentially distributed time.
In \cite{BoxmaResing}, the joint queue length process is analysed both at
embedded time points (beginnings of glue periods, visit periods and
switch-over periods) and at arbitrary time points, for the model with two
queues and {\it deterministic} glue periods. This analysis is later on
extended in Abidini, Boxma and Resing \cite{Abidini16} to the model with a
general number of queues. After that, in Abidini et al. \cite{KimKim}, an algorithm is presented to obtain the moments of
the number of customers in each station for the model with {\it
exponentially distributed} glue periods. Furthermore, in
\cite{KimKim} also
a workload decomposition for the model with {\it generally distributed}
glue periods is derived leading to a pseudo-conservation law.
The pseudo-conservation law in its turn is used to obtain approximations
of the mean waiting times at the different stations. In these papers, however, no analytical expressions for the complete joint distributions have been derived, which is something we aim to do in this paper.


In this manuscript, we will study the above-described polling system with
glue periods and retrials in a heavy traffic regime. More concretely, we will regard the regime where each of the arrival rates is scaled with the same constant, and subsequently the constant approaches from below that value, for which the system is critically loaded. Then, the workload
offered to the server is scaled to such a proportion that the queues are on
the verge of instability. Many techniques have been used to obtain the
heavy traffic behaviour of a variety of different polling models. Initial
studies of the heavy traffic behaviour of polling systems can be found in
Coffman, Puhalskii and Reiman \cite{CoffmanEtAl1,CoffmanEtAl2}, where the
occurrence of a so-called heavy traffic averaging principle is established.
This principle implies that, although the total scaled workload in the system
tends to a Bessel-type diffusion in the heavy traffic regime, it may be
considered as a constant during the course of a polling cycle, while the loads
of the individual queues fluctuate like a fluid model. It will turn out that
this principle will also hold true for this model. Furthermore, in van der Mei
\cite{MeiHTLST}, several heavy traffic limits have been established by taking
limits in known expressions for the Laplace-Stieltjes transform (LST) of the
waiting-time distribution. Alternatively, Olsen and van der Mei \cite{OlsenMei}
provide similar results, by studying the behaviour of the descendant set approach
(a numerical computation method, cf.\ Konheim, Levy and Srinivasan
\cite{KonheimLevySrinivasan}) in the heavy traffic limit. For the derivation of
heavy traffic asymptotics for our model, however, we will use results from
branching theory, mainly those presented in Quine \cite{Quine}.  Earlier, these
results have resulted in heavy traffic asymptotics for conventional polling models,
see van der Mei \cite{MeiHTBranching}. We will use the same method as presented
in that paper, but for a different class of polling system that models the dynamics
of optical networks. In addition, for some steps of the analysis, we will present new and straight forward proofs, while other steps require a different approach. Furthermore, we will derive asymptotics for the \emph{joint} queue length process at arbitrary time points, as opposed to just the marginal processes as derived in \cite{MeiHTBranching}. Due to the additional intricacies of the model at hand, we
will need to overcome many arising complex difficulties, as will become apparent later.

The rest of the paper is organized as follows. In Section 2, we introduce
some notation and present a theorem from \cite{Quine} on multitype
branching processes with immigration. In Section 3, we describe in detail
the polling model with retrials and glue periods and recall from
\cite{Abidini16} how the joint queue length process at
some embedded time points in this model is related to multitype branching
processes with immigration. Next, we will derive heavy traffic results for
our model. In Section 4, we consider the joint queue length process at
the start of glue periods. In Section 5, we look at the joint queue
length process at the start of visit and switch-over periods, while in
Section 6, we consider the joint queue length
process at arbitrary time points. Finally, in Section 7, we show how the heavy
traffic results, in combination with a light traffic result, can be used
to approximate performance measures for stable systems with arbitrary
system loads.

\section{Multitype branching processes with immigration}

To derive heavy-traffic results for the model under study, we regard its queue length process as a multitype branching process with immigration. To this end, before introducing the actual model in detail, we will state an important general result from \cite{Quine} on multitype branching process with immigration in this section, which we will make significant use of in the sequel of this paper. To state this result, we will first need some notation.

A multitype branching process with immigration has two kinds of individuals: immigrants and offspring.
We denote with $\underbar{$g $} = (g_1, \ldots, g_n)$ the mean immigration
vector. Here, $g_i$ is the mean number of type~$i$ immigrants entering
the system in each generation, for all $i= 1,\ldots, N$.
The offspring in the model is represented by the vector of generating functions
$h(\underbar{$z$}) = ( h_{1}(\underbar{$z$}), h_{2}(\underbar{$z$}), \ldots, h_{N}(\underbar{$z$}))$.
Here, $\underbar{$z$}=(z_1,z_2,\ldots,z_N)$ and  $|z_i| \leq 1$, for all
$i= 1, \ldots, N$, and
\begin{equation*}
h_i(\underbar{$z$}) = \sum_{j_1,\ldots,j_N \geq 0} p_i(j_1,\ldots,j_N) z_1^{j_1} \ldots z_N^{j_N},
\end{equation*}
where $p_i(j_1,\ldots,j_N)$ is the probability that a type~$i$ individual produces $j_k$ type~$k$ individuals, for all $i =1,\ldots,N$ and $k =1,\ldots,N$.
We use this to define the mean matrix $\mathbf{M} = (m_{i,j})$, where $m_{i,j}= \left.\frac{\partial h_i(\underbar{$z$})}{\partial z_j}\right\vert_{\underbar{$z$} = \underbar{$1$}}$, for all $i,j= 1, \ldots, N$, where \underbar{$1$} represents a vector of which each of the entries equals one. The elements $m_{ij}$ represent
the mean number of type $j$ children produced by a type~$i$ individual per generation.
We also define the second-order derivative matrix $K^{(i)}= \left(k^{(i)}_{j,k}\right)$ where $k^{(i)}_{j,k}= \left.\frac{\partial^2 h_i(\underbar{$z$})}{\partial z_j \partial z_k}\right\vert_{\underbar{$z$}=\underbar{$1$}}$, for all $i,j,k= 1, \ldots, N.$

Define  $\underbar{$ w$}  = (w_1, \ldots, w_N)^T$ as the normalized right eigenvector corresponding to the maximal
eigenvalue $\xi$ of $\mathbf{M}$. Then,
\[
 \mathbf{M} \underbar{$ w$}  =  \xi \underbar{$w $} ~~~~~~~~~~~~ \text{and}~~~~~~~~~~~~~\underbar{$ w$}  ^T \underbar{$1 $}  = 1.
\]

Furthermore, we define $\underbar{$v$} = (v_1, \ldots, v_N)^T$  as the left eigenvector of $\mathbf{M}$, corresponding to the maximal
eigenvalue $\xi$, normalized such that
\[
 \underbar{$v $} ^T \underbar{$ w$} =1.
\]

Additionally we give the following general notation in order to state the result of \cite{Quine}. Any variable $x$ which is dependent on $\xi$ will be denoted by $\hat{x}$ to indicate that it is
evaluated at $\xi =1$. Further, for $0<\xi <1$ let
\begin{equation}
 \pi_0(\xi) :=0~~~~~ \text{and}~~~~~ \pi_n(\xi) := \sum_{r=1}^{n} \xi^{r-2}, ~~~~~~n = 1,2,\ldots.
\label{lifeofpi}
\end{equation}

We denote with $\Gamma (\alpha, \mu)$ a gamma-distributed random variable. For  $\alpha , \mu , x > 0$, its probability density function is given by
\[
 f(x) = \frac{\mu^{\alpha}}{\Gamma(\alpha)} x^{\alpha -1} e^{- \mu x},~~~~~
 \text{where} ~~~~~~ \Gamma(\alpha) := \int_{t=0}^{\infty} t^{\alpha -1} e^{-t} dt.
\]

Now that the required notation is defined, we state the following important result, which is Theorem 4 of \cite{Quine}.

\begin{theorem}
If all first and second order derivatives of $h(\underbar{$z$})$  exist at $\underbar{$ z$} = \underbar{$ 1$} $, and $ 0 <g_i < \infty$ for all $i = 1, \ldots, N$, then
 \begin{equation*}
  \frac{1}{\pi(\xi)}
   \left(\begin{array}{cc}
         Z_{1} \\
         \vdots \\
         Z_{N}
        \end{array} \right)
\xrightarrow[d]{}
A   \left(\begin{array}{cc}
         \hat{v}_1 \\
         \vdots \\
         \hat{v}_N
        \end{array} \right)
\Gamma(\alpha , 1),~~~\text{when}~~ \xi \uparrow 1 .
        \end{equation*}
Here $\xrightarrow[d]{}$ means convergence in distribution, $ \pi(\xi) :=\lim_{n \to \infty}\pi_n(\xi)$, $\alpha := \frac{1}{A}  \hat{\underbar{$ g$} }^T \hat{\underbar{$w $} }$ and $A:= \frac{1}{2} \sum_{i=1}^{N} \hat{v}_i \left( \hat{\underbar{$w $}} ^T \hat{K}^{(i)} \hat{\underbar{$ w$} }\right) > 0 $.
The vector $(Z_{1},Z_{2},\ldots ,Z_{N})$ is defined such that $Z_{i}$ is the steady-state number of individuals of type~$i$ in the mulitype branching process with immigration, for all $i= 1, \ldots, N.$
\label{thm1}
\end{theorem}

\section{Polling model with retrials and glue periods}
In this section we first define the polling model with retrials and glue
periods. Then, we recall from \cite{Abidini16} its property that the
joint queue length process at the start of glue periods of a certain queue
is a multitype branching process with immigration.

\subsection{Model description}
\label{subsec:model}
We consider a single server polling model
with multiple queues,
$Q_i,$ $i= 1,\ldots,N$.
Customers arrive at $Q_i$ according to a Poisson process with rate $\lambda_i$;
they are called type-$i$ customers.
The service times at $Q_i$ are i.i.d., with $B_i$ denoting a generic service time,
with distribution $B_i(\cdot)$ and Laplace-Stieltjes transform (LST) $\tilde{B}_i(\cdot)$.
The server cyclically visits all the queues, thus after a visit of $Q_i$,
it switches to $Q_{i+1}$, $i= 1,\ldots,N$.
Successive switch-over times from $Q_i$ to $Q_{i+1}$
are i.i.d., where $S_i$ denotes a generic switch-over time,
with distribution $S_i(\cdot)$ and LST $\tilde{S}_i(\cdot)$.
We make all the usual independence assumptions about interarrival times, service times
and switch-over times at the queues.
After a switch of the server to $Q_i$, there first is a deterministic (i.e., constant)
glue period $G_i$, before the visit of the server at $Q_i$ begins.
The significance of the glue period stems from the following assumption.
Customers who arrive at $Q_i$ do not receive service immediately.
When customers arrive at $Q_i$ during a glue period $G_i$, they stick, joining the queue of $Q_i$.
When they arrive in any other period, they immediately leave and
retry after a retrial interval which is independent
of everything else, and which is exponentially distributed with parameter $\nu_i$, $i=1,\ldots,N$.

Since customers will only `stick' during the glue period, the service discipline at all queues can be interpreted as being gated. That is, during the visit period at $Q_i$, the server serves all
`glued' customers in that queue, i.e., all type-$i$ customers waiting at the end of the glue period,
but none of those in orbit,
and neither any new arrivals. We are interested in the steady-state behaviour of this polling model with retrials.
We hence assume that the stability condition
$\rho = \sum_{i=1}^N \rho_i < 1$ holds, where
$\rho_i := \lambda_i \E [B_i]$.

Note that now the server has three different periods at each station, a deterministic glue period during which customers are glued for service, followed by a visit period during which all the glued
customers are served and a switch-over period during which the server moves to the next station.
We denote, for $i=1,\ldots ,N$, by $(X_{1}^{(i)},X_{2}^{(i)},\ldots ,X_{N}^{(i)})$, $(Y_1^{(i)},Y_2^{(i)},\ldots ,Y_N^{(i)})$ and
$(Z_1^{(i)},Z_2^{(i)},\ldots ,Z_N^{(i)})$ vectors with as distribution the limiting distribution of
the number of customers of the different types in the system at the start of a glue period, a visit period and a switch-over period of station $i$, respectively. Furthermore, we denote, for $i=1,\ldots ,N$, by $(V_1^{(i)},V_2^{(i)},\ldots ,V_N^{(i)})$
the vector with as distribution the limiting distribution of
the number of customers of the different types in the system at an arbitrary point in time during a visit period of station $i$.
During glue and visit periods, we furthermore distinguish between those customers who are queueing in $Q_i$
and those who are in orbit for $Q_i$. Therefore we write
$Y_i^{(i)} = Y_i^{(iq)} + Y_i^{(io)}$ and $V_i^{(i)} = V_i^{(iq)} + V_i^{(io)}$, for all $i=1,\ldots ,N$, where $q$ represents queueing and $o$ represents in orbit.
Finally we denote by $(L^{(1q)},\ldots ,L^{(Nq)},L^{(1o)},\ldots,L^{(No)})$ the vector with as distribution the limiting distribution of
the number of customers of the different types in the queue and in the orbit at an arbitrary point in time.

The generating function of the vector of numbers of arrivals
at $Q_1$ to $Q_N$ during a type-$i$ service time $B_i$ is $\beta_i(\underbar{$z$}) := \tilde{B}_i(\sum_{j=1}^{N}\lambda_j(1-z_j))$.
Similarly, the generating function of the vector of numbers of arrivals
at $Q_1$ to $Q_N$ during a type-$i$ switch-over time $S_i$ is $\sigma_i(\underbar{$z$}) := \tilde{S}_i(\sum_{j=1}^{N}\lambda_j(1-z_j))$.

\subsection{Relation with multitype branching processes}
\label{sub:branchingprocess}
We now identify the relation of the polling model as defined in Section \ref{subsec:model} with a multitype branching process.
In \cite{Abidini16} it is shown that the number of customers of different
types in the system at the start of a glue period of station $1$ in the
polling model with retrials and glue periods is a multitype branching
process with immigration. In particular, it is derived in \cite{Abidini16}
that the joint PGF of $X_1^{(1)},\ldots,X_N^{(1)}$ satisfies
\begin{equation}
 \E \left[z_1^{X_{1}^{(1)}} z_2^{X_{2}^{(1)}}\ldots z_N^{X_{N}^{(1)}}\right]=\prod_{i=1}^{N}\sigma^{(i)}(\underbar{$z$}) \prod_{i=1}^{N} {\rm e}^{-G_i D_i(\underbar{$z$})}
 \E \left[ [h_1(\underbar{$z$})]^{X_{1}^{(1)}}[h_2(\underbar{$z$})]^{X_{2}^{(1)}}\ldots[h_N(\underbar{$z$})]^{X_{N}^{(1)}}\right],
 \label{M1}
 \end{equation}
 \begin{eqnarray}
 \text{where}~~~~~~\sigma^{(i)}(\underbar{$z$}) &:=& \sigma_i(z_1,\ldots ,z_i,h_{i+1}(\underbar{$z$}),\ldots ,h_N(\underbar{$z$})), \nonumber \\
 D_i(\underbar{$z$})&:=& \sum_{j=1}^{i-1} \lambda_j (1-z_j) +\lambda_i \left(1-\beta^{(i)}(\underbar{$z$})\right) + \sum_{j=i+1}^{N} \lambda_j(1-h_j(\underbar{$z$})), \nonumber \\
\beta^{(i)}(\underbar{$z$}) &:=& \beta_i(z_1,\ldots ,z_i,h_{i+1}(\underbar{$z$}),\ldots ,h_N(\underbar{$z$})), \nonumber \\
 h_{i}(\underbar{$z$})&:=&f_i(z_1,\ldots ,z_i,h_{i+1}(\underbar{$z$}),\ldots ,h_N(\underbar{$z$})),\nonumber \\
   \text{and}~~~~~~f_i(\underbar{$z$}) &:=& (1-{\rm e}^{-\nu_i G_i}) \beta_i(\underbar{$z$}) + {\rm e}^{-\nu_i G_i} z_i. \nonumber
 \end{eqnarray}

As explained in detail in \cite{Abidini16}, \eqref{M1} consists of the
product of three factors:
\begin{itemize}
\item $\prod_{i=1}^{N}\sigma^{(i)}(\underbar{$z$})$ represents new arrivals during
switch-over times and descendants of these arrivals in the current cycle.
\item  $\prod_{i=1}^{N} {\rm e}^{-G_i D_i(\underbar{$z$})}$ represents new arrivals during glue periods
and descendants of these arrivals in the current cycle.
The function $D_i(\underbar{$z$})$ is itself a sum of three terms:
\begin{itemize}
 \item $\sum_{j=1}^{i-1} \lambda_j (1-z_j)$ represents the arrivals of type $j <i$; these arrivals are not served in the current cycle.
 \item $\lambda_i \left(1-\beta^{(i)}(\underbar{$z$})\right)$ represents descendants of the arrivals of type $i$; these arrivals are all served
during the visit of station $i$ in the current cycle.
\item  $\sum_{j=i+1}^{N} \lambda_j(1-h_j(\underbar{$z$}))$ represents the
arrivals or descendants of arrivals of type $j>i$; these arrivals are either served (with probability $1-{\rm e}^{-\nu_i G_i}$) or not served (with probability
${\rm e}^{-\nu_i G_i}$) in the current cycle.
\end{itemize}
\item $\E \left[ [h_1(\underbar{$z$})]^{X_{1}^{(1)}}[h_2(\underbar{$z$})]^{X_{2}^{(1)}}\ldots[h_N(\underbar{$z$})]^{X_{N}^{(1)}}\right]$ represents descendants of $(X_1^{(1)},\ldots,X_N^{(1)})$ generated in the current cycle.
 \end{itemize}

We now proceed to further identify the branching process by finding its mean matrix $M$ and the mean immigration vector \underbar{$g$}.

\subsubsection*{Mean matrix of branching process:}

The elements $m_{i,j}$ of the mean matrix $\mathbf{M}$ of the branching
process are given by
\begin{equation*}
m_{i,j}= f_{i,j} \cdot 1[j \leq i] + \sum_{k=i+1}^{N} f_{i,k} m_{k,j},
\label{meanchildre}
\end{equation*}
where $f_{i,j}= \left.\frac{\partial f_i(\underbar{$z$})}{\partial z_j}\right\vert_{\underbar{$z$}=\underbar{$1$}},$ and hence
\begin{eqnarray*}
 f_{i,j} & = & \begin{cases}
          (1- {\rm e}^{-\nu_iG_i}) \lambda_j \E [B_i] , &  i \ne j, \\
       (1- {\rm e}^{-\nu_iG_i}) \rho_i + {\rm e}^{-\nu_iG_i}  , & i=j.
      \end{cases}
    \end{eqnarray*}

In the heavy traffic analysis of this model the following lemma will be useful.

 \begin{lemma}
\begin{equation}
 \mathbf{M} = \mathbf{M_1} \cdots \mathbf{M_N},
 \label{piM}
\end{equation}
where, for $i=1,2,\ldots,N$, we have
 \begin{equation}
\mathbf{M_i} = \left(\begin{array}{ccccccccc}
1~~~&~~~ 0~~~&~~~ \cdots~~~&~~~ 0 ~~~&~~~0 ~~~&~~~0~~~&~~~ \cdots ~~~&~~~\cdots~~~&~~~ 0 \\
0 ~~~&~~~1 ~~~&~~~\ddots~~~&~~~ \vdots~~~&~~~ 0~~~&~~~ 0 ~~~&~~~\cdots~~~&~~~ \cdots ~~~&~~~0 \\
\vdots~~~&~~~\ddots~~~&~~~\ddots~~~&~~~ 0~~~&~~~ 0~~~&~~~ 0~~~&~~~ \cdots~~~&~~~ \cdots~~~&~~~ 0\\
0~~~&~~~\cdots~~~&~~~0~~~&~~~ 1~~~&~~~ 0~~~&~~~ 0~~~&~~~ \cdots~~~&~~~ \cdots~~~&~~~ 0\\
 f_{i,1}& f_{i,2}&\cdots& f_{i, i-1}& f_{i,i}& f_{i, i+1}& \vdots& \vdots& f_{i,N}\\
0~~~&~~~ \cdots ~~~&~~~\cdots~~~&~~~ 0  ~~~&~~~ 0 ~~~&~~~ 1 ~~~&~~~ 0~~~&~~~ \cdots ~~~&~~~ 0    \\
0~~~&~~~ \cdots ~~~&~~~\cdots~~~&~~~ 0  ~~~&~~~ 0 ~~~&~~~ 0 ~~~&~~~ 1   ~~~&~~~ \ddots ~~~&~~~ 0 \\
0~~~&~~~ \cdots ~~~&~~~\cdots~~~&~~~ 0  ~~~&~~~ 0 ~~~&~~~ 0   ~~~&~~~ \ddots ~~~&~~~ \ddots ~~~&~~~ 0   \\
0~~~&~~~ \cdots ~~~&~~~\cdots~~~&~~~ 0   ~~~&~~~ 0 ~~~&~~~ 0   ~~~&~~~ \cdots ~~~&~~~ 0 ~~~&~~~ 1  \\
 \end{array} \right),
\label{meanmatrixi}
\end{equation}

\end{lemma}
\proof{First of all note that $ m_{N,j}= f_{N,j}$ for all $j= 1, \ldots, N.$ Therefore we have

\begin{equation*}
\mathbf{M_N} = \left(\begin{array}{ccccccccc}
1~~~&~~~ 0~~~&~~~ \cdots~~~&~~~\cdots~~~&~~~ 0 \\
0 ~~~&~~~1 ~~~&~~~\ddots~~~&~~~ \cdots ~~~&~~~0 \\
\vdots~~~&~~~\ddots~~~&~~~\ddots~~~&~~~ \cdots~~~&~~~ 0\\
0~~~&~~~ \cdots ~~~&~~~\cdots~~~&~~~ \ddots ~~~&~~~ 0   \\
0~~~&~~~ \cdots ~~~&~~~\cdots~~~&~~~ 1 ~~~&~~~0  \\
 m_{N,1}& m_{N,2}&\cdots&  m_{N,N-1}& m_{N,N}\\
 \end{array} \right).
\end{equation*}
Now using the fact that
$m_{N-1,j}= f_{N-1,j} + f_{N-1,N} m_{N,j}$ for all $j \leq N-1$ and
furthermore
$m_{N-1,N}= f_{N-1,N} m_{N,N}$, we obtain that
\begin{eqnarray*}
\mathbf{M_{N-1}}\mathbf{M_N}
 &=& \left(\begin{array}{ccccccccc}
1~~~&~~~ 0~~~&~~~ \cdots~~~&~~~ \cdots ~~~&~~~\cdots~~~&~~~ 0 \\
0 ~~~&~~~1 ~~~&~~~\ddots~~~&~~~\cdots~~~&~~~ \cdots ~~~&~~~0 \\
\vdots~~~&~~~\ddots~~~&~~~\ddots~~~&~~~ \cdots~~~&~~~ \cdots~~~&~~~ 0\\
0~~~&~~~ \cdots ~~~&~~~\cdots~~~&~~~ \ddots ~~~&~~~ \ddots ~~~&~~~ 0   \\
0~~~&~~~ \cdots ~~~&~~~\cdots~~~&~~~ 1 ~~~&~~~ 0 ~~~&~~~0  \\
 m_{N-1,1}& m_{N-1,2}&\cdots& \cdots& m_{N-1,N-1}&  m_{N-1,N}\\
 m_{N,1}& m_{N,2}&\cdots& \cdots& m_{N,N-1}& m_{N,N}\\
 \end{array} \right).
\end{eqnarray*}
Continuing in this way we obtain
\begin{eqnarray*}
\mathbf{M_{1}} \cdots \mathbf{M_N} = \left(\begin{array}{cc}
m_{1,1} \cdots m_{1,N} \\
\vdots ~~~\ddots ~~~\vdots\\
m_{N,1} \cdots m_{N,N}
\end{array} \right) = \mathbf{M}.
\end{eqnarray*}
\qed

\begin{remark}
(Intuition behind Lemma 1) The matrix $\mathbf{M_i}$ represents what happens with customers during a visit period at station~$i$. Customers at station~$i$ itself are either served
or not served, leading to the $i^{th}$ row with elements $f_{i,j}$. Customers at all other stations are not served leading to $1$'s on the diagonal and $0$'s
outside the diagonal. We obtain the product $\mathbf{M_{1}} \cdots \mathbf{M_N}$ because a cycle consists successively of visit periods of station~$1$, station~$2$, $\ldots$, up to station~$N$.
\end{remark}

\subsubsection*{Mean number of immigrants:}

Next, we look at the immigration part of the process. Let $g_i$ be the
mean number of type~$i$ individuals which immigrate into the system in
each generation. Equation (3.12) of \cite{Abidini16} gives us

\begin{equation}
g_i =\sum_{k=1}^{N}  \lambda_k  \Bigg(\Bigg(\sum_{j=1}^{k-1}(G_j+\E [S_j])\Bigg)\big(1- {\rm e}^{-\nu_k G_k}\big) + G_k\Bigg) m_{k,i}+\lambda_i \left(\sum_{j=1}^{i-1}(G_j+\E [S_j]){\rm e}^{-\nu_i G_i}+ \sum_{j=i}^N \E [S_j] + \sum_{j=i+1}^N G_j\right) .
\label{immigrants}
\end{equation}

The right-hand side of (\ref{immigrants}) is the sum of two terms. The term $\sum_{k=1}^{N}  \lambda_k  \left(\left(\sum_{j=1}^{k-1}(G_j+\E [S_j])\right)\left(1- {\rm e}^{-\nu_k G_k}\right) + G_k\right) m_{k,i}$ represents the mean number of type~$i$
customers which are descendants of customers of type $k$, arriving during glue periods and switch-over periods before the visit period of station $k$ and served during the visit at station~$k$, in the current cycle.
The first part of the second term $\lambda_i \sum_{j=1}^{i-1}(G_j+\E [S_j]){\rm e}^{-\nu_i G_i}$ represents the mean number of customers of type~$i$ which
arrive during glue periods and switch-over periods before the visit of the server at station~$i$ and which are not served during the visit of station $i$ in the current cycle.
The second part of the second term $\lambda_i \left(\sum_{j=i}^N \E [S_j] + \sum_{j=i+1}^N G_j\right)$ represents the mean number of customers of type~$i$ which
arrive during glue periods and switch-over periods after the visit period of station~$i$ in the current cycle.
Note that each of the terms mentioned above is non-negative and finite.
Furthermore, for non-zero glue periods and arrival rates, at least one of
the terms is non-zero. Therefore we have $0< g_i < \infty$.

\begin{remark}
Note that the branching part of the process only represents descendants of customers which are present in the system at the start of a glue period of station $1$.
Customers which arrive at stations during glue periods and switch-over periods are not represented by the branching part of the process. Instead they and
their descendants are represented by the immigration
part of the process. Both glue periods and switch-over periods can be
considered as parts of the cycle during which the server is not working. This rather unexpected feature explains why the polling model at hand is not part of the class of polling models considered by \cite{MeiHTBranching}, but requires an analysis on its own.
\end{remark}

\section{Heavy traffic analysis: number of customers at start of glue periods of station $1$}

Now that we have successfully modelled the polling system as a multitype branching process with immigration, we derive the limiting scaled joint queue length distribution in each station at the start of glue periods of station $1$ by following the same line of proof as that of \cite{MeiHTBranching}.
In \cite{MeiHTBranching}, the author first proves a couple of lemmas for a conventional branching-type polling system without retrials and glue periods and, afterwards, uses
these lemmas to give the heavy traffic asymptotics of the joint queue length process at certain embedded time points. In the following subsection, we will derive similar lemmas in order to derive a heavy traffic theorem for our polling system with retrials and glue periods.

Note that when we scale our system such that $\rho \uparrow 1$, we
are effectively changing the arrival rate at each station while keeping
the service times and the ratios of the arrival rates fixed. Let any
variable $x$ which is dependent on $\rho$ be denoted by $\hat{x}$
whenever it is evaluated at $\rho =1$. Therefore we have for any system
that, $\lambda_i= \rho \hat{\lambda}_i $.

From Theorem 1 of \cite{Zedek}, we know that if all elements of a matrix
are continuous in some variable, then the real eigenvalues of this matrix
are also continuous in that variable. As each element of $\mathbf{M}$ is a
continuous function of $\rho$, the maximal eigenvalue $\xi$
is a continuous function of $\rho$ as well. Furthermore, from Lemmas 3, 4 and 5 of
\cite{Resing93} we know that $\xi<1$ when $\rho<1$, $\xi=1$ when $\rho =1$
and $\xi > 1$ when $\rho>1$. Therefore, we have that $\xi$ is a
continuous function of $\rho$ and
\begin{equation*}
  \lim_{\rho \uparrow 1}\xi (\rho) = \xi (1)=1.
  \label{xi1}
 \end{equation*}

\subsection{Preliminary results and lemmas}\label{sec:prelim}

\begin{lemma}
The normalized right and left eigenvectors of $\mathbf{\hat{M}}$, the mean matrix of the system with $\rho=1$, corresponding to the maximal
eigenvalue $\xi =1$, are respectively given by
 \begin{equation*}
  \hat{\underbar{$w $} } =  \left(\begin{array}{cc}
         \hat{w}_1 \\
         \vdots \\
         \hat{w}_N
        \end{array} \right) =\frac{ \underbar{$b $}  }{\lvert \underbar{$b $}  \rvert }~~~~~~~~~~~~~~~~~~~~~~\text{and}~~~~~~~~~~~~~~~~~~ \hat{\underbar{$ v$}}
     =  \left(\begin{array}{cc}
         \hat{v}_1 \\
         \vdots \\
         \hat{v}_N
        \end{array} \right) = \frac{\lvert \underbar{$ b$} \rvert} {\delta }  \hat{\underbar{$ u$}},
        \end{equation*}
  where
\begin{equation*}
   \underbar{$b $}  =  \left(\begin{array}{cc}
         \E [B_1] \\
         \vdots \\
         \E [B_N]
        \end{array} \right),~~~~~~~~~~~~~~~~~~~~~~~~~\underbar{$ u$}
     =  \left(\begin{array}{cc}
         u_1 \\
         \vdots \\
         u_N
        \end{array} \right),
        \end{equation*}

\[ \lvert \underbar{$ b$} \rvert : = \sum_{j=1}^N \E [B_j],~~~~~~~~~~u_j := \lambda_j \left[\frac{e^{-\nu_j G_j} }{1-e^{-\nu_j G_j}} + \sum_{k=j}^N \rho_k \right]~~~\text{and}~~~~\delta := \hat{\underbar{$ u$} }^{T}\underbar{$b $} .
\]
\label{lemma:eigenvector}
\end{lemma}
\proof{
First we look at the normalized right eigenvector $\hat{ \underbar{$w $} }$. Using \eqref{meanmatrixi} we evaluate the vector $  \mathbf{\hat{M}_i} \hat{\underbar{$ w$}  }$.
Let $  \left( \mathbf{\hat{M}_i} \hat{\underbar{$ w$}} \right)_j  $ represent the $j^{th}$ element of $  \mathbf{\hat{M}_i} \hat{\underbar{$ w$}  }$. By a series of simple algebraic manipulations, it follows then that
 \begin{equation*}
  \left( \mathbf{\hat{M}_i} \hat{\underbar{$ w$}} \right)_j  =
  \begin{cases} \hat{w}_j  , &  j \ne i, \\
      \frac{1}{\lvert \underbar{$ b$}  \rvert} \sum_{k=1}^{N}  \hat{f}_{i,k} \E[B_k] , & j=i.
      \end{cases}
 \end{equation*}
However, it also holds that
\begin{eqnarray*}
 \sum_{k=1}^{N} \hat{f}_{i,k}  \E[B_k] &=& e^{-\nu_i G_i} \E[B_i]+ \sum_{k=1}^{N} (1-e^{-\nu_i G_i}) \E [B_i] \hat{\lambda}_k  \E [B_k] \\
 &=&  e^{-\nu_i G_i} \E[B_i]+ (1-e^{-\nu_i G_i}) \E [B_i] \sum_{k=1}^{N} \hat{\rho}_k = \E [B_i] .
\end{eqnarray*}
Therefore, we conclude that $\left( \mathbf{\hat{M}_i} \hat{\underbar{$ w$}} \right)_i = \hat{w}_i$.
This implies that $\hat{ \underbar{$ w$}  }$ is the normalized right eigenvector of  $  \mathbf{\hat{M}_i}$ for an eigenvalue $\xi=1$, for all $i= 1, \ldots,N$.
Hence from \eqref{piM} we get the first part of the lemma.
Next, we look at the left eigenvector $\hat{\underbar{$ v$} }$. Since $\hat{\underbar{$u $}}$ is a multiple of $\hat{\underbar{$v $}}$, it is enough to show that $\hat{\underbar{$ u$} } $ is an eigenvector of $\hat{\mathbf{ M}} $.
Define
\begin{eqnarray}
\underbar{$ u$}^{(i)}
     =  \left(\begin{array}{cc}
         u^{(i)}_1 \\
         \vdots \\
         u^{(i)}_N
        \end{array} \right),~~~~~~~ \text{where}~~~~~~~~~~u_j^{(i)} =\begin{cases}
       \lambda_j \left[\frac{e^{-\nu_j G_j} }{1-e^{-\nu_j G_j}} + \sum_{k=j}^N \rho_k  + \sum_{k=1}^{i-1} \rho_k \right] , &  i \leq j ,\\
       \lambda_j \left[\frac{e^{-\nu_j G_j} }{1-e^{-\nu_j G_j}}  + \sum_{k=j}^{i-1} \rho_k \right] , &  i>j.
   \end{cases}
\label{eq:hat}
\end{eqnarray}
Note that $u^{(1)}_j = u^{(N+1)}_j = u_j$, for all $j= 1, \ldots, N$, and
hence, $\underbar{$u$}^{(1)}=\underbar{$u$}^{(N+1)}=\underbar{$u$}$.
Furthermore, we have
\begin{eqnarray*}
 \hat{\underbar{$ u$} }^{{(1)}^T} \mathbf{\hat{M}_1} =  \left(\begin{array}{cc}
                                                                \hat{u}_1 \hat{f}_{1,1}\\
                                                                 \hat{u}_1 \hat{f}_{1,2} +  \hat{u}_2\\
                                                                 \vdots\\
                                                                 \hat{u}_1 \hat{f}_{1,N} +  \hat{u}_N\\
                                                              \end{array}\right)^T =  \left(\begin{array}{cc}
                                                                \hat{\lambda}_1 \frac{e^{-\nu_1 G_1}}{1-e^{-\nu_1 G_1}} + \hat{\lambda}_1 \hat{\rho}_1\\
                                                                   \hat{u}_2  +\hat{\lambda}_2 \hat{\rho}_1\\
                                                                 \vdots\\
                                                                \hat{u}_N  + \hat{\lambda}_N \hat{\rho}_1\\
                                                              \end{array}\right)^T   =   \left(\begin{array}{cc}
                                                                \hat{u}_1^{(2)} \\
                                                                  \hat{u}_2^{(2)}\\
                                                                 \vdots\\
                                                                \hat{u}_N^{(2)}\\
                                                              \end{array}\right)^T = \hat{\underbar{$ u$} }^{{(2)}^T}     ,
\end{eqnarray*}
and, in a similar way, for all $i= 1, \ldots, N$,
\begin{equation}
\hat{\underbar{$ u$} }^{{(i)}^T} \mathbf{\hat{M}_i} = \hat{\underbar{$ u$} }^{{(i+1)}^T}.
\label{lefteigeni}
\end{equation}
Therefore we have
\[
\hat{\underbar{$u $}}^{T}  \mathbf{M} =\hat{\underbar{$ u$} }^{{(1)}^T} \mathbf{\hat{M}_1}\cdots \mathbf{\hat{M}_N} =
        \hat{\underbar{$ u$} }^{{(N+1)}^T} = \hat{\underbar{$ u$} }^{T}.
\]
Hence $\hat{\underbar{$ u$} }$ and  $\hat{\underbar{$ v$}}$ are the left eigenvectors of  $\mathbf{\hat{M}}$, for eigenvalue $\xi = 1$.
 \qed
\begin{remark}
Alternatively, we could have used lemma 4 from \cite{MeiHTBranching} to find the left
and normalized right eigenvectors. The normalized right eigenvector $\hat{\underbar{$w$}}$ is the same as given in \cite{MeiHTBranching}.
To find the left eigenvector $\hat{\underbar {$v$}}$ from \cite{MeiHTBranching}, we
first need to calculate the exhaustiveness factor $f_j$. In our model, this
exhaustiveness factor is given by
$f_j = (1- e^{-\nu_j G_j}) (1 -\rho_j)$.
Each customer of type~$j$, present at the start of a glue period at station~$j$, is served with probability $(1- e^{-\nu_j G_j})$ and during that
service time on average $\rho_j$ new type~$j$ customers will arrive.
Furthermore, with probability $e^{-\nu_j G_j}$ a customer of type~$j$,
present at the start of a glue period at station~$j$, is not served.
Therefore we have $1- f_j =  (1- e^{-\nu_j G_j}) \rho_j+ e^{-\nu_j G_j} $, and hence the exhaustiveness factor is given by $f_j = (1- e^{-\nu_j G_j}) (1 -\rho_j).$
Substituting this exhaustiveness factor in lemma 4 of \cite{MeiHTBranching} we get
\begin{eqnarray*}
 u_j &=& \lambda_j \left[\frac{(1-\rho_j)(1-(1- e^{-\nu_j G_i}) (1 -\rho_j)) }{(1- e^{-\nu_j G_j}) (1 -\rho_j)} + \sum_{k=j+1}^N \rho_k \right] \nonumber \\
 &=& \lambda_j \left[\frac{ e^{-\nu_j G_j} + (1 -e^{-\nu_j G_j})\rho_j }{1- e^{-\nu_j G_j} } + \sum_{k=j+1}^N \rho_k \right] \nonumber \\
 &=& \lambda_j \left[\frac{ e^{-\nu_j G_j} }{1- e^{-\nu_j G_j} } + \sum_{k=j}^N \rho_k \right],
\end{eqnarray*}
which is in agreement with Lemma \ref{lemma:eigenvector}.

\end{remark}

\begin{remark}
\label{rem:vector}
 In Lemma \ref{lemma:eigenvector} we have given the left and normalized right eigenvectors for the mean matrix $\mathbf{\hat{M}}$ at eigenvalue $\xi =1$.
Note that this mean matrix is defined for the branching process when we
consider the beginning of a glue period of station~$1$ as
the initial point of the cycle.
Instead if we consider
the beginning of a glue period of station~$i$ as the initial point of the
cycle, we get, for eigenvalue $\xi =1$, the same normalized right
eigenvector $\hat{\underbar{$w$}}$. However, the left eigenvector is now given by the vector $\hat{\underbar{$ v$}}^{(i)}$ defined by
  \[
   \hat{\underbar{$ v$}}^{(i)}
     =  \left(\begin{array}{cc}
         \hat{v}^{(i)}_1 \\
         \vdots \\
         \hat{v}^{(i)}_N
        \end{array} \right) = \frac{\lvert \underbar{$ b$} \rvert} {\delta }  \hat{\underbar{$ u$}}^{(i)}.
    \]
    Note that $\delta = \hat{\underbar{$ u$} }^{(1)^T}\underbar{$b $}  = \hat{\underbar{$ u$} }^{(i)^T}\underbar{$b $}.$
In this paper we prove all the lemmas and theorems using $ \hat{\underbar{$ v     $}}= \hat{\underbar{$ v$}}^{(1)}$. However, we can instead use $ \hat{\underbar{$ v$}}^{(i)}$ and prove the same by just changing the initial
    point of the cycle from the beginning of a glue period of station~$1$ to the beginning of a glue period of station~$i$.
\end{remark}

In Lemma 2 we have evaluated the normalized right, and left eigenvectors of $\mathbf{M}$ at the maximal eigenvalue, when $\rho \uparrow 1$.
We will now use this to compute the value of the derivative of this eigenvalue as $\rho \uparrow 1$.

\begin{lemma}
For the maximal eigenvalue $\xi = \xi (\rho)$ of the matrix $\mathbf{M}$, the derivative of $\xi (\rho)$ w.r.t. $\rho$ satisfies
\begin{equation*}
 \xi ^ \prime (1)  = \frac{1}{\delta}.
\label{xiderivative}
 \end{equation*}
\label{lemma:eigenvaluederi}
\end{lemma}
\proof{
Since the maximal eigenvalue $\xi$ of $\mathbf{M}$ is a simple eigenvalue
and furthermore $\mathbf{M}$ is continuous in $\rho$, Theorem 5 of Lancaster \cite{Lancaster64} states that
\begin{equation}
 \left.\frac{d \xi}{d \rho} \right\vert_ {\rho=1} = \frac{  \hat{\underbar{$ v$} }^{T} \mathbf{\hat{M}^{\prime}}\hat{\underbar{$w $}} }{  \hat{\underbar{$v $} }^{T}\hat{\underbar{$ w$}}},
 \label{derivative}
\end{equation}
where $\mathbf{\hat{M}^{\prime}}$ is the element wise derivative of $\mathbf{M}$ w.r.t. $\rho$ evaluated at $\rho =1$.
Let $\mathbf{U_i} =\left( \prod_{k=1}^{i-1}  \mathbf{M_k} \right) \mathbf{M^{\prime}_i} \left(\prod_{k=i+1}^{N}  \mathbf{M_k}\right)$. Then due to \eqref{piM} we can write
  $\mathbf{M^{\prime}} = \sum_{i=1}^{N}  \mathbf{U_i}$.
From \eqref{meanmatrixi} we can see that
 \begin{eqnarray}
  \mathbf{M^{\prime}_i} &=& \left(\begin{array}{cccccc}
0~~~& \cdots~~~&~~~ 0 ~~~&~~~\cdots~~~&~~~ 0 \\
\vdots~~~&~~~\ddots~~~&~~~ \cdots~~~&~~~ \ddots~~~&~~~ \vdots\\
 \frac{df_{i,1}}{d\rho}&\cdots& \frac{df_{i,i}}{d\rho}& \cdots& \frac{df_{i,N}}{d\rho}\\
\vdots~~~&~~~ \ddots ~~~&~~~\cdots~~~&~~~ \vdots   \\
0~~~&~~~ \cdots ~~~&~~~ 0   ~~~&~~~ \cdots ~~~&~~~ 0  \\
 \end{array} \right)\nonumber \\
 &=& \left(\begin{array}{cccc}
 0~~~& ~~~\cdots~~~&~~~ 0 \\
 \vdots~~~&~~~ \ddots~~~&~~~ \vdots\\
(1- {\rm e}^{-\nu_iG_i}) \E[B_i]\frac{d\lambda_1}{d\rho}&~~~ \cdots~~~&(1- {\rm e}^{-\nu_iG_i}) \E[B_i] \frac{d\lambda_N}{d\rho} \\
\vdots~~~&~~~ \ddots~~~&~~~ \vdots\\
 0~~~&~~~ \cdots~~~&~~~0  \\
 \end{array} \right)\nonumber\\
 &=&  \left(\begin{array}{cc}
~~~~0 \\
~~~~\vdots\\
(1- {\rm e}^{-\nu_iG_i}) \E[B_i]\\
~~~~\vdots\\
~~~~0 \\
 \end{array} \right) \left(   \frac{d\lambda_1}{d\rho}~~~\cdots~~~  \frac{d\lambda_N}{d\rho} \right).
  \label{midder}
 \end{eqnarray}
From the definition of $\rho$, we know that
$\sum_{i=1}^N \E [B_i]  \frac{d\lambda_i}{d\rho} = 1$, and hence
\begin{equation}
\left(\frac{d \lambda_1}{d\rho}~~~\cdots~~~  \frac{d\lambda_N}{d\rho} \right) \hat{\underbar{$w $}}= \lvert \underbar{$b $} \rvert ^{-1}.
\label{finalder}
\end{equation}
Since $\hat{\underbar{$w $}}$ is the normalized right eigenvector of any $ \mathbf{\hat{M}_i}$ for eigenvalue $\xi=1$, we have
\begin{equation}
 \prod_{k=i+1}^{N}  \mathbf{\hat{M}_k} \hat{\underbar{$w $}} =  \hat{\underbar{$w $}}.
 \label{initialder}
\end{equation}
Using \eqref{midder}, \eqref{finalder} and \eqref{initialder} we get
\begin{equation}
  \mathbf{\hat{M}^{\prime}_i} \prod_{k=i+1}^{N}  \mathbf{\hat{M}_k} \hat{\underbar{$ w$} } = \frac{1}{\lvert \underbar{$ b$}  \rvert}   \left(\begin{array}{cc}
~~~~0 \\
~~~~\vdots\\
(1- {\rm e}^{-\nu_iG_i}) \E[B_i]\\
~~~~\vdots\\
~~~~0 \\
 \end{array} \right).
 \label{lemma2p1}
\end{equation}
From \eqref{lefteigeni} and \eqref{lemma2p1} we get
\begin{eqnarray}
 \hat{\underbar{$u $}}^{T} \mathbf{\hat{U}_i} \hat{\underbar{$w $} } &=&
 \left(\hat{\underbar{$ u$} }^{T} \prod_{k=1}^{i-1}  \mathbf{\hat{M}_k}  \right)  \mathbf{\hat{M}^{\prime}_i} \left(\prod_{k=i+1}^{N}  \mathbf{\hat{M}_k} \hat{\underbar{$ w$} } \right) \nonumber \\
 &=& \frac{1}{\lvert \underbar{$ b$}  \rvert}    \left(\begin{array}{cc}
      \hat{u}_{1}^{(i)} \\
      \vdots \\
      \hat{u}_{i-1}^{(i)} \\
           \hat{u}_i^{(i)} \\
        \hat{u}_{i+1}^{(i)} \\
                \vdots \\
          \hat{u}_N^{(i)}
        \end{array} \right) ^T \left(\begin{array}{cc}
~~~~0 \\
~~~~\vdots\\
~~~~0 \\
(1- {\rm e}^{-\nu_iG_i}) \E[B_i]\\
~~~~0 \\
~~~~\vdots\\
~~~~0 \\
 \end{array} \right) \nonumber \\
 &=& \frac{\hat{u}_{i}^{(i)} (1- {\rm e}^{-\nu_iG_i}) \E[B_i]}{\lvert \underbar{$ b$}  \rvert} = \frac{\hat{\rho}_i}{\lvert \underbar{$ b$}  \rvert},
 \label{lemma2p2}
\end{eqnarray}
where the last equality follows from the fact that
$\hat{u}_{i}^{(i)} = \hat{\lambda}_i/(1- {\rm e}^{-\nu_iG_i})$, see \eqref{eq:hat}.
Multiplying both sides of \eqref{lemma2p2} with $\lvert \underbar{$ b$} \rvert/ \delta$ and summing it over all $i= 1,\ldots,N,$ we get that
\begin{equation}
\hat{\underbar{$ v$}}^{T} \mathbf{\hat{M}^\prime} \hat{\underbar{$w $} } =
\sum_{i=1}^{N} \frac{\lvert \underbar{$ b$} \rvert }{\delta}\hat{\underbar{$u $} }^{T} \mathbf{\hat{U}_i} \hat{\underbar{$w $} } = \sum_{i=1}^{N} \frac{\hat{\rho}_i}{\delta} = \frac{1}{\delta}.
\label{eq:lemma2}
\end{equation}
Since $\hat{\underbar{$ v$}}^{T}  \hat{\underbar{$w $} } = 1$, we obtain from \eqref{derivative} and \eqref{eq:lemma2} that
$ \xi ^ \prime (1) =  \frac{1}{\delta}.$
\qed}

\vspace{0.2cm}

In Theorem~\ref{thm1}, we need all the second-order derivatives
$\frac{\delta^2 h_{i}(\underbar{$z$})}{\delta z_j \delta z_k}$
of the function $h_{i}(\underbar{$z$})$.
In Lemma 4, we first find $\frac{\delta^2 h_{i}(\underbar{$z$})}{\delta z_j \delta z_k}$, for all $i,j$ and $k$, and
then use them to find the parameter $A$ as defined in Theorem~\ref{thm1}.

\begin{lemma}
 For the second-order derivative matrix $K^{(i)}= \left(k^{(i)}_{j,k}\right)$ where $k^{(i)}_{j,k}= \left. \frac{\partial^2 h_i(\underbar{$z$})}{\partial z_j \partial z_k}\right\vert_{\underbar{$z$} =\underbar{$1$}}$, for all $i,j,k= 1, \ldots, N,$
 we have that \[
A:= \frac{1}{2}\sum_{i=1}^{N} \hat{v}_i^{(1)} \left( \hat{\underbar{$w $}} ^T \hat{K}^{(i)} \hat{\underbar{$ w$} }\right) =  \frac{1}{2 \delta \lvert \underbar{$b $}  \rvert}\frac{b^{(2)}}{b^{(1)}},
                 \]
                 where $b^{(j)} = \frac{\sum_{i=1}^{N} \lambda_i \E [B_i^j]}{\sum_{i=1}^{N} \lambda_i},$ for $j=1,2$.
\label{lemma:Avalue}
                 \end{lemma}
\proof{

We know that
\begin{eqnarray*}
  h_{i}(\underbar{$z$})&=&f_i(z_1,\ldots ,z_i,h_{i+1}(\underbar{$z$}),\ldots ,h_N(\underbar{$z$})) \\
&=&  (1-{\rm e}^{-\nu_i G_i}) \beta_i(z_1,\ldots ,z_i,h_{i+1}(\underbar{$z$}),\ldots ,h_N(\underbar{$z$})) + {\rm e}^{-\nu_i G_i} z_i\\
&=& (1-{\rm e}^{-\nu_i G_i}) \E \left[ {\rm e}^{-B_i\left(\sum_{c=1 }^{i}(1-z_c) \lambda_c  + \sum_{c=i+1}^{N}(1-h_c(\underbar{$z$})) \lambda_c \right)}\right]+ {\rm e}^{-\nu_i G_i} z_i.\end{eqnarray*}
From this it follows that
\begin{eqnarray*}
\frac{\partial h_{i}(\underbar{$z$})}{\partial z_k} &=& (1-{\rm e}^{-\nu_i G_i}) \E \left[B_i\left( \lambda_k  1[k \leq i] + \sum_{c=i+1}^{N}  \lambda_c \frac{\partial h_{c}(\underbar{$z$})}{\partial z_k} \right) {\rm e}^{-B_i\left(\sum_{c=1 }^{i}(1-z_c) \lambda_c  + \sum_{c=i+1}^{N}(1-h_c(\underbar{$z$})) \lambda_c \right)}\right]+ {\rm e}^{-\nu_i G_i} 1[k = i] ,\\
\text{and}~~~~ \frac{\partial^2 h_{i}(\underbar{$z$})}{\partial z_j \partial z_k} &=& (1-{\rm e}^{-\nu_i G_i}) \E \Bigg[\Bigg(B_i^2 \left( \lambda_k  1[k \leq i] + \sum_{c=i+1}^{N}  \lambda_c \frac{\partial h_{c}(\underbar{$z$})}{\partial z_k} \right)\left( \lambda_j  1[j \leq i] + \sum_{c=i+1}^{N}  \lambda_c \frac{\partial h_{c}(\underbar{$z$})}{\partial z_j} \right)\\
&&~~~~~~~~~~~~~~~~~~~~+B_i \sum_{c=i+1}^{N}  \lambda_c \frac{\partial^2 h_{c}(\underbar{$z$})}{\partial z_j \partial z_k} \Bigg) {\rm e}^{-B_i\left(\sum_{c=1 }^{i}(1-z_c) \lambda_c  + \sum_{c=i+1}^{N}(1-h_c(\underbar{$z$})) \lambda_c \right)}\Bigg],
\end{eqnarray*}
where $1[E]= 1$, when the event $E$ is true and otherwise $1[E]= 0$.

Because $\left. \frac{\partial^2 h_i(\underbar{$z$})}{\partial z_j \partial z_k}\right\vert_{\underbar{$z$} =\underbar{$1$}}=k^{(i)}_{j,k} $ and $(1-{\rm e}^{-\nu_i G_i})\E[B_i]\left( \lambda_k  1[k \leq i] + \sum_{c=i+1}^{N}  \lambda_c \frac{\delta h_{c}}{\delta z_k}(i\underbar{$1$})\right)  = m_{i,k} - 1[k = i]{\rm e}^{-\nu_i G_i}$, we have
\begin{eqnarray}
k^{(i)}_{j,k} &=& \frac{\E [B_i^2]}{\E[B_i]^2(1-{\rm e}^{-\nu_i G_i})}(m_{i,j}- 1[j = i]{\rm e}^{-\nu_i G_i}) (m_{i,k}- 1[k = i]{\rm e}^{-\nu_i G_i}) + (1-{\rm e}^{-\nu_i G_i}) \E[B_i] \sum_{c=i+1}^{N} \lambda_c k^{(c)}_{j,k} \nonumber \\
&=&\frac{\E [B_i^2]}{\E[B_i]^2(1-{\rm e}^{-\nu_i G_i})} (m_{i,j} m_{i,k} -(1[j = i] m_{i,k} +1[k = i]m_{i,j}){\rm e}^{-\nu_i G_i} +1[i= j =k ]{\rm e}^{-2\nu_i G_i} )\nonumber\\
&&+ (1-{\rm e}^{-\nu_i G_i}) \E[B_i] \sum_{c=i+1}^{N} \lambda_c k^{(c)}_{j,k}. 
\label{eq:kijk}
\end{eqnarray}
Let $\mathbf{1}_i$ be an $N\times N$ matrix, where the element in the $i$-th row and the $i$-th column equals one, and all $N^2-1$ other entries read zero. Then, based on \eqref{eq:kijk}, we can write
\begin{eqnarray}
 K^{(i)} &=& \frac{\E [B_i^2]}{\E[B_i]^2(1-{\rm e}^{-\nu_i G_i})} \left[\left(\begin{array}{cc}
m_{i,1}m_{i,1} \cdots m_{i,1}m_{i,N} \\
\vdots ~~~\ddots ~~~\vdots\\
m_{i,N}m_{i,1} \cdots m_{i,N} m_{i,N}
\end{array} \right)  -{\rm e}^{-\nu_i G_i} \left(\begin{array}{cc}
0 ~~\cdots~~ m_{i,1} ~~\cdots~~0\\
\vdots ~~~\ddots ~~~\vdots~~~\ddots ~~~\vdots\\
0 ~\cdots ~m_{i,i-1} ~\cdots~ 0\\
m_{i,1}~ \cdots ~2m_{i,i}~\cdots~ m_{i,N}\\
0 ~\cdots ~m_{i,i+1} ~\cdots ~0\\
\vdots ~~~\ddots ~~~\vdots~~~\ddots ~~~\vdots\\
0 ~~\cdots~~ m_{i,N}~~\cdots~~ 0
\end{array} \right) +{\rm e}^{-2\nu_i G_i} \mathbf{1}_i
 \right] \nonumber\\
 &&+(1-{\rm e}^{-\nu_i G_i}) \E[B_i] \sum_{c=i+1}^{N} \lambda_cK^{(c)} \nonumber\\
 &=& \frac{\E [B_i^2]}{\E[B_i]^2(1-{\rm e}^{-\nu_i G_i})} \left[\left(\begin{array} {cc}
m_{i,1} \\  \vdots \\ m_{i,N} \end{array} \right)\left(\begin{array}{cc}
m_{i,1} \cdots m_{i,N} \end{array} \right)  - {\rm e}^{-\nu_i G_i}\left( \begin{array}{cc}
0 ~~\cdots~~ m_{i,1} ~~\cdots~~0\\
\vdots ~~~\ddots ~~~\vdots~~~\ddots ~~~\vdots\\
0 ~\cdots ~m_{i,i-1} ~\cdots~ 0\\
m_{i,1}~ \cdots ~2m_{i,i}~\cdots~ m_{i,N}\\
0 ~\cdots ~m_{i,i+1} ~\cdots ~0\\
\vdots ~~~\ddots ~~~\vdots~~~\ddots ~~~\vdots\\
0 ~~\cdots~~ m_{iN}~~\cdots~~ 0
\end{array} \right) + {\rm e}^{-2\nu_i G_i} \mathbf{1}_i
 \right] \nonumber\\
 &&+(1-{\rm e}^{-\nu_i G_i}) \E[B_i] \sum_{c=i+1}^{N} \lambda_c K^{(c)}. \nonumber \end{eqnarray}
This leads to\begin{eqnarray}
\underbar{$w$}^T K^{(i)} \underbar{$w$} &=&\frac{\E [B_i^2]}{\E[B_i]^2(1-{\rm e}^{-\nu_i G_i})} \left[\underbar{$w$}^T \left(\begin{array} {cc}
m_{i,1} \\  \vdots \\ m_{i,N} \end{array} \right)\left(\begin{array}{cc}
m_{i,1} \cdots m_{i,N} \end{array} \right) \underbar{$w$}  - {\rm e}^{-\nu_i G_i}\underbar{$w$}^T \left( \begin{array}{cc}
0 ~~\cdots~~ m_{i,1} ~~\cdots~~0\\
\vdots ~~~\ddots ~~~\vdots~~~\ddots ~~~\vdots\\
0 ~\cdots ~m_{i,i-1} ~\cdots~ 0\\
m_{i,1}~ \cdots ~2 m_{i,i}~\cdots~ m_{i,N}\\
0 ~\cdots ~m_{i,i+1} ~\cdots ~0\\
\vdots ~~~\ddots ~~~\vdots~~~\ddots ~~~\vdots\\
0 ~~\cdots~~ m_{i,N}~~\cdots~~ 0
\end{array} \right) \underbar{$w$} + {\rm e}^{-2\nu_i G_i} \underbar{$w$}^T \mathbf{1}_i \underbar{$w$} \right] \nonumber\\
 &&+(1-{\rm e}^{-\nu_i G_i}) \E[B_i] \sum_{c=i+1}^{N} \lambda_c \underbar{$w$}^T K^{(c)}\underbar{$w$}.
 \label{eq:lemma4_0}
\end{eqnarray}
Note that from the definition of $\hat{\underbar{$w$}}$, we have that
\begin{equation}
\hat{\underbar{$w$}}^T \left(\begin{array}{cc}
\hat{m}_{i,1} \\  \vdots \\ \hat{m}_{i,N} \end{array} \right)
= \left(\begin{array}{cc}
\hat{m}_{i,1} \cdots \hat{m}_{i,N} \end{array} \right) \hat{\underbar{$w$}}
= \frac{\E[B_i]}{\lvert \underbar{$b $}  \rvert } .
\label{eq:lemma4_1}
\end{equation}
Now we evaluate \begin{eqnarray}
               \hat{  \underbar{$w$}}^T \left( \begin{array}{cc}
0 ~~\cdots~~  \hat{m}_{i,1} ~~\cdots~~0\\
\vdots ~~~\ddots ~~~\vdots~~~\ddots ~~~\vdots\\
0 ~\cdots ~ \hat{m}_{i,i-1} ~\cdots~ 0\\
 \hat{m}_{i,1}~ \cdots ~2  \hat{m}_{i,i}~\cdots~  \hat{m}_{i,N}\\
0 ~\cdots ~ \hat{m}_{i,i+1} ~\cdots ~0\\
\vdots ~~~\ddots ~~~\vdots~~~\ddots ~~~\vdots\\
0 ~~\cdots~~  \hat{m}_{i,N}~~\cdots~~ 0
\end{array} \right) \hat{\underbar{$w$}} &=& \frac{1}{\lvert \underbar{$b $}  \rvert } \left( \begin{array}{cc}  \hat{m}_{i,1}\E[B_i] \\ \vdots  \\\hat{m}_{i,i-1}\E[B_i]\\ \hat{m}_{i,i}\E[B_i] + \sum_{j=1}^{N} \hat{m}_{i,j} \E[B_j]\\ \hat{m}_{i,i+1}\E[B_i]\\ \vdots \\ \hat{m}_{i,N}\E[B_i] \end{array}\right)^T\hat{\underbar{$w$}} \nonumber \\
             &=&  \frac{\E[B_i]\sum_{j=1}^{N} \hat{m}_{i,j}\E[B_j] + \E[B_i]\sum_{j=1}^{N} \hat{m}_{i,j}\E[B_j]}{\lvert \underbar{$b $}  \rvert ^2} \nonumber \\
             &=& \frac{2 \E[B_i]^2} {\lvert \underbar{$b $}  \rvert ^2} .
             \label{eq:lemma4_2}
             \end{eqnarray}
   Evaluating \eqref{eq:lemma4_0} for $\rho \uparrow 1$, and substituting \eqref{eq:lemma4_1} and \eqref{eq:lemma4_2} in it, we have
   \begin{eqnarray}
    \hat{\underbar{$w$}}^T \hat{K}^{(i)} \hat{\underbar{$w$}} &=&\frac{\E [B_i^2]}{\E[B_i]^2(1-{\rm e}^{-\nu_i G_i})} \left[\frac{\E[B_i]^2}{\lvert \underbar{$b $}  \rvert ^2}- 2 {\rm e}^{-\nu_i G_i}\frac{\E[B_i]^2}{\lvert \underbar{$b $}  \rvert ^2} + {\rm e}^{-2\nu_i G_i}\frac{\E[B_i]^2}{\lvert \underbar{$b $}  \rvert ^2} \right] +(1-{\rm e}^{-\nu_i G_i}) \E[B_i] \sum_{c=i+1}^{N} \hat{\lambda}_c \hat{\underbar{$w$}}^T \hat{K}^{(c)} \hat{\underbar{$w$}} \nonumber \\
   &=& (1-{\rm e}^{-\nu_i G_i})\frac{\E [B_i^2]}{\lvert \underbar{$b $}  \rvert ^2} +(1-{\rm e}^{-\nu_i G_i}) \E[B_i] \sum_{c=i+1}^{N} \hat{\lambda}_c \hat{\underbar{$w$}}^T \hat{K}^{(c)} \hat{\underbar{$w$}} \nonumber \\
   &=&(1-{\rm e}^{-\nu_i G_i}) \left(\frac{\E [B_i^2]}{\lvert \underbar{$b $}  \rvert ^2} +\E[B_i]\sum_{c=i+1}^{N} \hat{\lambda}_c \hat{\underbar{$w$}}^T \hat{K}^{(c)} \hat{\underbar{$w$}}  \right).
   \label{lemma4_3}
   \end{eqnarray}
Multiplying both sides of \eqref{lemma4_3} with $\hat{v}_i$ and evaluating it for $i=1$ we get
\begin{eqnarray}
\hat{v}_1 \hat{\underbar{$w$}}^T \hat{K}^{(1)} \hat{\underbar{$w$}} &=&\frac{ \lvert \underbar{$ b$} \rvert \hat{\lambda}_1}{\delta} \left(\frac{\E [B_1^2]}{\lvert \underbar{$b $}  \rvert ^2} +\E[B_1]\sum_{c=2}^{N} \hat{\lambda}_c \hat{\underbar{$w$}}^T \hat{K}^{(c)} \hat{\underbar{$w$}}  \right)\nonumber \\
&=&\frac{ \hat{\lambda}_1 \E [B_1^2]}{\delta \lvert \underbar{$b $}  \rvert } + \frac{\lvert \underbar{$ b$} \rvert\hat{\rho}_1 \hat{\lambda}_2 (1-{\rm e}^{-\nu_2 G_2})}{\delta} \left(\frac{\E [B_2^2]}{\lvert \underbar{$b $}  \rvert ^2} +\sum_{c=3}^{N} \hat{\lambda}_c \hat{\underbar{$w$}}^T \hat{K}^{(c)} \hat{\underbar{$w$}}  \right) + \lvert \underbar{$ b$} \rvert\frac{\hat{\rho}_1}{\delta}\sum_{c=3}^{N} \hat{\lambda}_c \hat{\underbar{$w$}}^T \hat{K}^{(c)} \hat{\underbar{$w$}},
\label{eq:lemma4_4}
\end{eqnarray}
where for the second equality we again used \eqref{lemma4_3}, but now for $i=2$,
to substitute $\hat{\underbar{$w$}}^T \hat{K}^{(2)} \hat{\underbar{$w$}}$.
Multiplying both sides of \eqref{lemma4_3} with $\hat{v}_i$ and evaluating it for $i=2$ we get
\begin{eqnarray}
\hat{v}_2 \hat{\underbar{$w$}}^T \hat{K}^{(2)} \hat{\underbar{$w$}} &=& \frac{\lvert \underbar{$ b$} \rvert\hat{\lambda}_2}{\delta} \left( {\rm e}^{-\nu_2 G_2}+(1-{\rm e}^{-\nu_2 G_2}) \sum_{j=2}^{N} \hat{\rho}_j \right) \left(\frac{\E [B_2^2]}{\lvert \underbar{$b $}  \rvert ^2} +\E[B_2]\sum_{c=3}^{N} \hat{\lambda}_c \hat{\underbar{$w$}}^T \hat{K}^{(c)} \hat{\underbar{$w$}}  \right)\nonumber\\
&=& \frac{\lvert \underbar{$ b$} \rvert\hat{\lambda}_2}{\delta}\left( {\rm e}^{-\nu_2 G_2}+(1-{\rm e}^{-\nu_2 G_2}) (1-\hat{\rho}_1) \right) \left(\frac{\E [B_2^2]}{\lvert \underbar{$b $}  \rvert ^2} +\E[B_2]\sum_{c=3}^{N} \hat{\lambda}_c \hat{\underbar{$w$}}^T \hat{K}^{(c)} \hat{\underbar{$w$}}  \right)\nonumber\\
&=& \frac{\lvert \underbar{$ b$} \rvert\hat{\lambda}_2 -  \lvert \underbar{$ b$} \rvert\hat{\rho}_1\hat{\lambda}_2 (1-{\rm e}^{-\nu_2 G_2})}{\delta}  \left(\frac{\E [B_2^2]}{\lvert \underbar{$b $}  \rvert ^2} +\E[B_2]\sum_{c=3}^{N} \hat{\lambda}_c \hat{\underbar{$w$}}^T \hat{K}^{(c)} \hat{\underbar{$w$}}  \right).
\label{eq:lemma4_5}
\end{eqnarray}
Summing \eqref{eq:lemma4_4} and \eqref{eq:lemma4_5} we get
\begin{equation*}
\sum_{j=1}^{2} \hat{v}_j   \hat{\underbar{$w$}}^T \hat{K}^{(j)} \hat{\underbar{$w$}} = \sum_{j=1}^{2} \frac{ \hat{\lambda}_j}{\delta}\frac{\E [B_j^2]}{\lvert \underbar{$b $}  \rvert } +  \frac{\lvert \underbar{$ b$} \rvert}{\delta}\left(\sum_{j=1}^{2} \hat{\rho}_j \right)\sum_{c=3}^{N}\hat{\lambda}_c \hat{\underbar{$w$}}^T \hat{K}^{(c)} \hat{\underbar{$w$}}.
\end{equation*}
By repeating the above procedure, we end up with
\begin{equation*}
\sum_{j=1}^{N} \hat{v}_j \hat{\underbar{$w$}}^T \hat{K}^{(j)} \hat{\underbar{$w$}} = \sum_{j=1}^{N} \frac{ \hat{\lambda}_j}{\delta}\frac{\E [B_j^2]}{\lvert \underbar{$b $}  \rvert }=  \frac{1}{\delta \lvert \underbar{$b $}  \rvert} \frac{b^{(2)}}{b^{(1)}} .
\end{equation*}
Therefore we have \[
                   A:= \frac{1}{2}\sum_{j=1}^{N} \hat{v}_j \hat{\underbar{$w$}}^T \hat{K}^{(j)} \hat{\underbar{$w$}} = \frac{1}{2 \delta \lvert \underbar{$b $}  \rvert} \frac{b^{(2)}}{b^{(1)}} .
                  \]
\qed
}

At this point, we have determined all parameters required to deploy Theorem~\ref{thm1}, except for the constant $\alpha$. This parameter depends on the immigration part of our process and is given by the following lemma.

\begin{lemma}
 For $g = (g_1 \cdots g_N)^T$, we have that
 \begin{equation}
 \alpha:= \frac{1}{A} \hat{\underbar{$ g$} }^T \hat{\underbar{$w $} } = 2 r \delta \frac{b^{(1)} }{b^{(2)}},
\label{normalizedimmigrants}
 \end{equation}
where $r= \sum_{i=1}^{N} \left( \E[S_i] +G_i \right)$.
 \label{lemma:alpha}
\end{lemma}
\proof {

Multiplying both sides of \eqref{immigrants} with $\E [B_i]$ and summing it over all $i$ gives
\begin{eqnarray*}
 \sum _{i=1}^{N} g_i \E [B_i] &=& \sum_{k=1}^{N}  \lambda_k  \left(\sum_{i=1}^{N}  m_{k,i} \E [B_i]\right) \left(\sum_{j=1}^{k-1}\left(G_j+\E [S_j]\right)(1- {\rm e}^{-\nu_k G_k}) + G_k\right) \\
 &+& \sum_{i=1}^{N} \rho_i \left(\sum_{j=1}^{i-1}\left(G_j+\E [S_j]\right){\rm e}^{-\nu_i G_i}+ \sum_{j=i}^N \E [S_j] + \sum_{j=i+1}^N G_j\right).
 \end{eqnarray*}
Since $\hat{w}$ is an eigenvector of $\mathbf{\hat{M}_k}$, we have $\sum_{i=1}^{N} \hat{m}_{k,i} \E [B_i] = \E [B_k]$.
Hence, taking $\rho \uparrow 1$, we get

 \begin{eqnarray}
 \sum_{i=1}^{N} \hat{g_i} \E [B_i] &=& \sum_{i=1}^{N}  \hat{\rho_i}  \left(\sum_{j=1}^{i-1}\left(G_j+\E [S_j]\right)(1- {\rm e}^{-\nu_i G_i}) + G_i
 +  \sum_{j=1}^{i-1}(G_j+\E [S_j]){\rm e}^{-\nu_i G_i}+ \sum_{j=i}^N \E [S_j] + \sum_{j=i+1}^N G_j\right) \nonumber \\
 &=& \sum_{i=1}^{N}  \hat{\rho_i} \sum_{j=1}^{N} \left(\E [S_j] + G_j \right) = \sum_{i=1}^{N} (\E [S_i] + G_i).
 \label{eq:immi}
 \end{eqnarray}
Substituting $\hat{\underbar{$ w$} } = \lvert \underbar{$ b$}  \rvert ^{-1} (\E [B_1] ~\cdots ~ \E [B_N])^T$ and $A=   \frac{1}{2 \delta \lvert \underbar{$b $}  \rvert} \frac{b^{(2)}}{b^{(1)}}$ in \eqref{normalizedimmigrants}
and using  \eqref{eq:immi} will give that
$\alpha=  2 r \delta \frac{b^{(1)}}{b^{(2)}}$. \qed

\subsection{The heavy traffic theorem}

Similar to the procedure used in \cite{MeiHTBranching}, we will now combine the preliminary work in Section \ref{sec:prelim} with Theorem \ref{thm1} in order to obtain the following heavy traffic theorem for the complete queue length process at cycle starts.

\begin{theorem}\label{thm:startglue}
 For the cyclic polling system with retrials and glue periods, the scaled steady-state joint queue length vector at the start of glue periods at station~$1$ satisfies

 \begin{equation}
  (1-\rho)    \left(\begin{array}{cc}
         X_1^{(1)} \\
         \vdots \\
         X_N^{(1)}
        \end{array} \right)
\xrightarrow[d]{\rho \uparrow 1}\frac{b^{(2)}}{2b^{(1)}}\frac{1}{\delta}   \left(\begin{array}{cc}
         \hat{u}^{(1)}_{1} \\
         \vdots \\
         \hat{u}^{(1)}_{N}
        \end{array} \right)
\Gamma(\alpha , 1),
 \end{equation}
 where $\alpha = 2  r \delta\frac{b^{(1)}}{b^{(2)}}.$
 \label{thm:2}
\end{theorem}
\proof{
The joint queue length process at the start of glue periods of station~$1$ is an $N-$dimensional multitype branching process with mean matrix $\mathbf{M}$
and the mean number of type~$i$ immigration customers per generation, $g_i$, given by \eqref{piM} and \eqref{immigrants} respectively. At the end of Section \ref{sub:branchingprocess} we concluded that $0<g_i<\infty$.
Furthermore, $h(\underbar{$z$})= ( h_{1}(\underbar{$z$}), h_{2}(\underbar{$z$}), \ldots, h_{N}(\underbar{$z$}))$ is the offspring function for which the second-order derivatives $k^{(i)}_{j,k}= \left. \frac{\partial^2 h_i(\underbar{$z$})}{\partial z_j \partial z_k}\right\vert_{\underbar{$z$} =\underbar{$1$}}$, for all $i,j,k = 1,\ldots,N$, exist.
Therefore, all the conditions of Theorem \ref{thm1} are satisfied.
\\
Denote by $(X_{1}^{(1)},X_{2}^{(1)},\ldots ,X_{N}^{(1)})$
the vector with as distribution the limiting distribution of
the number of customers of the different types in the system at
the start of a glue period of station $1$.
Now, from Theorem \ref{thm1} it follows that
 \begin{equation}
  \frac{1}{\pi(\xi(\rho))}
   \left(\begin{array}{cc}
         X^{(1)}_{1} \\
         \vdots \\
         X^{(1)}_{N}
        \end{array} \right)
\xrightarrow[d]{}A   \left(\begin{array}{cc}
         \hat{v}_1^{(1)} \\
         \vdots \\
         \hat{v}_N^{(1)}
        \end{array} \right)
\Gamma(\alpha , 1), ~~~\text{when}~~ \rho \uparrow 1 ,
\label{thm2:1}
        \end{equation}
where $\pi(\xi(\rho)):=\lim_{n \to \infty}\pi_n(\xi(\rho))$,
and $A$ and $\hat{\underbar{$v$}}^{(1)}$ and
$\alpha = \frac{1}{A} \hat{\underbar{$ g$} }^T \hat{\underbar{$w $} },$
are as defined in
Lemmas \ref{lemma:eigenvector}, \ref{lemma:Avalue} and \ref{lemma:alpha}.

From \eqref{lifeofpi} we can say that, for $\rho < 1$,
\begin{equation*}
 \pi(\xi(\rho)) = \frac{1}{\xi(\rho) (1-\xi(\rho))}.
\end{equation*}
Using this, together with Lemma \ref{lemma:eigenvaluederi}, gives
\begin{equation}
 \lim_{\rho \uparrow 1} (1-\rho) \pi(\xi(\rho)) =  \lim_{\rho \uparrow 1} \frac{1 - \rho}{\xi(\rho) (1-\xi(\rho))} = \lim_{\rho \uparrow 1} \frac{-1 }{\xi^{\prime}(\rho) (1-2\xi(\rho))} = \lim_{\rho \uparrow 1} \frac{1}{\xi^{\prime}(\rho)}=\delta.
 \label{eq:limeps}
\end{equation}
Therefore, multiplying and dividing the LHS of \eqref{thm2:1} with $1-\rho$, we get
\begin{eqnarray*}
  \frac{1- \rho}{(1-\rho)\pi(\xi(\rho))}
   \left(\begin{array}{cc}
         X^{(1)}_{1} \\
         \vdots \\
         X^{(1)}_{N}
        \end{array} \right)
\xrightarrow[d]{}& A   \left(\begin{array}{cc}
         \hat{v}_1^{(1)} \\
         \vdots \\
         \hat{v}_N^{(1)}
        \end{array} \right)
\Gamma(\alpha , 1),~~~\text{when}~~ \rho \uparrow 1.
\end{eqnarray*}
Using \eqref{eq:limeps}, this gives
\begin{eqnarray*}
(1- \rho)
   \left(\begin{array}{cc}
         X^{(1)}_{1} \\
         \vdots \\
         X^{(1)}_{N}
        \end{array} \right)
\xrightarrow[d]{}&  \frac{1}{2 \lvert \underbar{$b $}  \rvert} \frac{b^{(2)}}{ b^{(1)}}  \left(\begin{array}{cc}
         \hat{v}_1^{(1)} \\
         \vdots \\
         \hat{v}_N^{(1)}
        \end{array} \right)
\Gamma(\alpha , 1), ~~~\text{when}~~ \rho \uparrow 1,
\end{eqnarray*}
and hence
\begin{eqnarray*}(1- \rho)
   \left(\begin{array}{cc}
         X^{(1)}_{1} \\
         \vdots \\
         X^{(1)}_{N}
        \end{array} \right)
\xrightarrow[d]{}&  \frac{b^{(2)}}{2b^{(1)}}\frac{1}{\delta}
\left(\begin{array}{cc}
         \hat{u}^{(1)}_1 \\
         \vdots \\
         \hat{u}^{(1)}_N
        \end{array} \right) \Gamma(\alpha , 1),~~~\text{when}~~ \rho \uparrow 1  .
 \end{eqnarray*}
\qed }



\subsection{Discussion of results: connection with a binomially gated polling model}
It turns out that the heavy traffic results that we obtained in this section for the model at hand, are similar to those of a binomially gated polling model (see e.g.\ \cite{Levy2}). The dynamics of the binomially gated polling model are much like those of a conventional gated polling model, except that after dropping a gate at $Q_i$, the customers before it will each be served in the corresponding visit period with probability $p_i$ in an i.i.d.\ way, rather than with probability one as in the gated model.
In particular, the heavy traffic analysis of our model coincides with that of a binomially gated polling model with the same
 interarrival time distributions, service time distributions and switch-over time distributions, and probability parameters $p_i=1-e^{-\nu_iG_i}$. To check this, we note that the binomially gated polling model with these probability parameters falls within the framework of the seminal work of \cite{MeiHTBranching} when taking the exhaustiveness parameters $f_i = (1-\rho_i)(1-e^{-\nu_iG_i})$, after which it is easily verified that Theorem 5 of \cite{MeiHTBranching} coincides with Theorem \ref{thm:startglue}. Note, however, that although we also exploit a branching framework in this paper, the model considered in this paper does not fall directly in the class of polling models considered in \cite{MeiHTBranching}, due to the intricate immigration dynamics it exposes.

The intuition behind this remarkable connection is as follows. First, we have that a binomially gated polling model does not have the feature of glue periods. However, in a heavy-traffic regime, the server in our model will reside in a visit period for 100\% of the time, so that glue periods hardly occur in this regime either. Furthermore, in a binomially gated polling model, each customer present at the start of a visit period at $Q_i$ will be served within that visit period with probability $p_i = 1-e^{-\nu_iG_i}$ in an i.i.d.\ fashion. Note that something similar happens with the model at hand. There, the start of a visit period coincides with the conclusion of a glue period. During this glue period, all customers present in the orbit of the queue will, independently from one another, queue up for the next visit period with probability $1-e^{-\nu_iG_i}$. These two facts explain the analogy.

Do note that this analogy, remarkable though it is, does not help us in the further analysis towards the asymptotics of the customer population at an arbitrary point in time. While Theorem \ref{thm:startglue} is now aligned with Theorem 5 of \cite{MeiHTBranching}, we cannot use the subsequent analysis steps in that paper to get to results concerning the customer population in heavy traffic at an arbitrary point in time. This is much due to the fact that the strategy of \cite{MeiHTBranching} exploits a relation between the queue length of $Q_1$ at a cycle start and the virtual waiting time of that queue at an arbitrary point in time. Since the type-$i$ customers in our model are not served in the order of arrival, as is usually assumed, such a relation is hard to derive and is essentially unknown. As an alternative, we will extend the current heavy traffic asymptotics at cycle starts to certain other embedded epochs in Section \ref{sec:embeddedTimePoints}, and eventually to arbitrary points in time in Section \ref{sec:arbitraryTimePoints}.

\section{Heavy traffic analysis: number of customers at other embedded time points}\label{sec:embeddedTimePoints}

A cycle in the polling system with retrials and glue periods passes
through three different phases: glue periods, visit periods and
switch-over periods. In the previous section, in Theorem \ref{thm:2},
we studied the behaviour of the scaled steady-state joint queue length
vector at the start of glue periods at station~$1$.
We will now extend this result to the scaled steady-state joint queue
length vector at the start of a visit period and the start of a
switch-over period in Theorems \ref{startofvisit} and \ref{startofswitch}.

\begin{theorem}
For the cyclic polling system with retrials and glue periods,
the scaled steady-state joint queue length vector
at the start of visit periods at station $1$ satisfies
\begin{equation}
  (1-\rho)    \left(\begin{array}{cc}
  Y_1^{(1q)}\\
        Y_1^{(1o)} \\
        Y_2^{(1)}\\
         \vdots \\
         Y_N^{(1)}
        \end{array} \right)
\xrightarrow[d]{}\frac{b^{(2)}}{2b^{(1)}}\frac{1}{\delta}   \left(\begin{array}{cc}
       (1-e^{-\nu_1 G_1})  \hat{u}^{(1)}_{1} \\
     e^{-\nu_1 G_1}  \hat{u}^{(1)}_{1} \\
       \hat{u}^{(1)}_{2} \\
        \vdots \\
         \hat{u}^{(1)}_{N}
        \end{array} \right)
\Gamma(\alpha , 1),~~~\text{when}~~ \rho \uparrow 1.
\label{eq:startvisit}
 \end{equation}
 \label{startofvisit}
\end{theorem}
\proof{

The distribution of the number of new customers of type~$j$ entering the system during a glue period of station~$i$ is stochastically smaller than that of the number of events $G_j^{(i)}$ in a Poisson process with rate $\hat{\lambda}_j$ during an interval of length $G_i$. This is due to the fact that the arrival rate $\lambda_j = \rho\hat{\lambda}_j$ does not exceed $\hat{\lambda}_j$.
Since $G_j^{(i)}$ is finite with probability 1, we have that
$(1-\rho)G_j^{(i)} \rightarrow 0$ with probability 1, as $\rho \uparrow 1$.
Therefore the limiting scaled joint queue length distribution, for all customers other than type~$i$,  is same at the start of a glue period
and at the start of a visit period of station~$i$.

Furthermore, the $X_i^{(i)}$ customers of type~$i$, present in the system at the start of a glue period of station~$i$, join the queue, independently of each other, with probability $1-e^{-\nu_i G_i}$ during the glue period.
Let $\{U_i, i \ge 0\}$ be a series of i.i.d. random variables where $U_k$ indicates whether the $k$-th customer joins the queue or stays in orbit, for all $k= 1, \ldots,X_i^{(i)}$. More specifically, $U_k =1$ if the customer joins the queue, with probability $1-e^{-\nu_i G_i}$, and $U_k=0$ if the customer stays in orbit, w.p. $e^{-\nu_i G_i}$. Then the number of customers of type~$i$ in the queue ($Y_i^{(iq)}$) and in the orbit
($Y_i^{(io)}$) at the start of a visit period at station~$i$ are given by
\begin{eqnarray*}
Y_i^{(iq)} = \sum_{k=1}^{ X_i^{(i)}} U_k  ~~~~~~&\text{and}&~~~~~~~ Y_i^{(io)} =  X_i^{(i)} -\sum_{k=1}^{ X_i^{(i)}} U_k.
\end{eqnarray*}
Since $X_i^{(i)} \to \infty$ with probability $1$, as $\rho \uparrow 1$, we have by virtue of the weak law of large numbers that
\begin{equation}
\frac{Y_i^{(iq)}}{X_i^{(i)}} = \frac{\sum_{k=1}^{ X_i^{(i)}} U_k}{ X_i^{(i)}}\xrightarrow[\mathbb{P}]{} 1- e^{-\nu_i G_i},~~~\text{when}~~ \rho \uparrow 1,
\label{eq:convergance1a}
\end{equation}
where $\xrightarrow[\mathbb{P}]{}$ means convergence in probability.
Similarly we have
\begin{equation}
\frac{Y_i^{(io)}}{X_i^{(i)}} = \frac{X_i^{(i)} -\sum_{k=1}^{ X_i^{(i)}} U_k}{ X_i^{(i)}}\xrightarrow[\mathbb{P}]{} e^{-\nu_i G_i},~~~\text{when}~~ \rho \uparrow 1.
\label{eq:convergance1b}
\end{equation}
Therefore, using Slutsky's convergence theorem \cite{Grimmett}, along with  \eqref{thm2:1}, \eqref{eq:convergance1a}, \eqref{eq:convergance1b} and the arguments above, we get
 \begin{equation*}
  (1-\rho)    \left(\begin{array}{cc}
  Y_1^{(1q)}\\
        Y_1^{(1o)} \\
        Y_2^{(1)}\\
         \vdots \\
         Y_N^{(1)}
        \end{array} \right)
\xrightarrow[d]{\rho \uparrow 1}\frac{b^{(2)}}{2b^{(1)}}\frac{1}{\delta}   \left(\begin{array}{cc}
       (1-e^{-\nu_1 G_1})  \hat{u}^{(1)}_{1} \\
     e^{-\nu_1 G_1}  \hat{u}^{(1)}_{1} \\
       \hat{u}^{(1)}_{2} \\
        \vdots \\
         \hat{u}^{(1)}_{N}
        \end{array} \right)
\Gamma(\alpha , 1), ~~~\text{when}~~ \rho \uparrow 1.
 \end{equation*}
\qed

We end this section by considering the scaled steady-state joint queue length vector at the start of a switch-over period from station~$1$ to station~$2$.

\begin{theorem}
For the cyclic polling system with retrials and glue periods,
the scaled steady-state joint queue length vector at the start of a
switch-over period from station~$1$ to station~$2$ satisfies
\begin{equation}
  (1-\rho)   \left(\begin{array}{cc}
         Z_1^{(1)} \\
         Z_2^{(1)} \\
         \vdots \\
         Z_N^{(1)}
        \end{array} \right)
\xrightarrow[d]{}\frac{b^{(2)}}{2b^{(1)}}\frac{1}{\delta}   \left(\begin{array}{cc}
        e^{-\nu_1 G_1} \hat{u}^{(1)}_{1} &+ (1-e^{-\nu_1 G_1}) \hat{u}^{(1)}_{1} \hat{\lambda_1} \E[B_1] \\
        \hat{u}^{(1)}_{2} &+ (1-e^{-\nu_1 G_1}) \hat{u}^{(1)}_{1} \hat{\lambda}_2 \E[B_1]\\
                 \vdots \\
         \hat{u}^{(1)}_{N} &+ (1-e^{-\nu_1 G_1}) \hat{u}^{(1)}_{1} \hat{\lambda}_N \E[B_1]
        \end{array} \right)
\Gamma(\alpha , 1), ~~~\text{when}~~ \rho \uparrow 1.
\label{eq:switch}
 \end{equation}
\label{startofswitch}
\end{theorem}

\proof{
The number of customers in the orbit of station~$j$ at the start of a switch-over period from station~$i$ to station~$i+1$ equals the number of customers in the orbit
at the start of the visit of station~$i$~ plus the Poisson arrivals with rate $\lambda_j$, during the service of customers in the queue of station~$i$, say $J_j^{(i)}$. In other words, we have that

\begin{equation}
 Z_j^{(i)} = \begin{cases}
        Y_j^{(i)} + J_j^{(i)}  , &  j \ne i, \\
       Y_i^{(io)} + J_i^{(i)} , & j=i.
      \end{cases}
      \label{startswitch1}
\end{equation}
Note that $J_j^{(i)}$ is the sum of Poisson arrivals with rate $\lambda_j$ during $Y_i^{(iq)}$ independent service times with distribution $B_i$. Let $D_{i,j,k}$ be the number of
Poisson arrivals with rate $\lambda_j$ during the $k^{th}$ service in the visit period of station $i$. Thus
\begin{equation*}
 J_j^{(i)} = \sum_{k=1}^{Y_i^{(iq)}} D_{i,j,k}.
\end{equation*}
Since $Y_i^{(iq)}\rightarrow \infty$ as $\rho \uparrow 1$, and $\E[B_i]$ is finite, we have by virtue of the weak law of large numbers that
\begin{equation}
 \frac{J_j^{(i)}}{ Y_i^{(iq)}} \xrightarrow[\mathbb{P}]{} \hat{\lambda}_j \E[B_i], ~~~\text{when}~~ \rho \uparrow 1.
 \label{startswitch2}
\end{equation}
Therefore, using Slutsky's convergence theorem along with \eqref{eq:startvisit}, \eqref{startswitch1} and \eqref{startswitch2} we get
\begin{equation*}
  (1-\rho)   \left(\begin{array}{cc}
         Z_1^{(1)} \\
         Z_2^{(1)} \\
         \vdots \\
         Z_N^{(1)}
        \end{array} \right)
\xrightarrow[d]{}\frac{b^{(2)}}{2b^{(1)}}\frac{1}{\delta}   \left(\begin{array}{cc}
        e^{-\nu_1 G_1} \hat{u}^{(1)}_{1} &+ (1-e^{-\nu_1 G_1}) \hat{u}^{(1)}_{1} \hat{\lambda_1} \E[B_1] \\
        \hat{u}^{(1)}_{2} &+ (1-e^{-\nu_1 G_1}) \hat{u}^{(1)}_{1} \hat{\lambda}_2 \E[B_1]\\
                 \vdots \\
         \hat{u}^{(1)}_{N} &+ (1-e^{-\nu_1 G_1}) \hat{u}^{(1)}_{1} \hat{\lambda}_N \E[B_1]
        \end{array} \right)
\Gamma(\alpha , 1), ~~~\text{when}~~ \rho \uparrow 1.
 \end{equation*}

\begin{remark}
Alternatively, Theorems \ref{startofvisit} and \ref{startofswitch} can be obtained by exploiting known relations between the joint PGFs of the vectors $(X_1^{(1)},\ldots, X_N^{(1)})$, $(Y_1^{(1q)},Y_1^{(1o)}, Y_2^{(1)}\ldots, Y_N^{(1)})$ and $(Z_1^{(1)},\ldots, Z_N^{(1)})$ given in Equations (3.2) and (3.3) of \cite{Abidini16}. After replacing each parameter $z_j$ in these functions by $z_j^{(1-\rho)}$ and taking the limit of $\rho$ going to one from below, these expressions give the relations between the joint PGFs of the heavy traffic distributions. Combining these results with Theorem \ref{thm:2} and subsequently invoking Levy's continuity theorem (see e.g.\ Section 18.1 of \cite{Williams}) then readily imply the theorems.
\end{remark}

\begin{remark}
Throughout this section, we have focused on the joint queue length process at the start of a glue,
visit or switch-over period at $Q_1$. However, similar results for the starts of these periods at any $Q_i$
can be obtained by either simply reordering indices, or by exploiting the relations obtained in \cite{Abidini16} between $(X_1^{(i)},\ldots, X_N^{(i)})$, $(Y_1^{(iq)},Y_1^{(io)}, Y_2^{(i)}\ldots, Y_N^{(i)})$, $(Z_1^{(i)},\ldots, Z_N^{(i)})$ and $(X_1^{(i+1)},\ldots, X_N^{(i+1)})$.
\end{remark}

\section{Heavy traffic analysis: number of customers at arbitrary time points}\label{sec:arbitraryTimePoints}

In this section we look at the limiting scaled joint queue length distribution of the number of customers at the different stations at an arbitrary time point. At such a point in time, the system
can be in the glue period, the visit period or the switch-over period of some station~$i$,
with probability $\frac{(1-\rho)G_i}{\sum_{i=1}^{N} (G_i+ \E[S_i])}$, $\rho_i$ and
$\frac{(1-\rho)\E[S_i]}{\sum_{i=1}^{N} (G_i+ \E[S_i])}$ respectively. As $\rho \uparrow 1$, the probabilities  $\frac{(1-\rho)G_i}{\sum_{i=1}^{N} (G_i+ \E[S_i])}$ and
$\frac{(1-\rho)\E[S_i]}{\sum_{i=1}^{N} (G_i+ \E[S_i])}$, both converge to $0$. Therefore we only need to study the scaled steady-state joint queue length vector at an arbitrary time
in each of the $N$ visit periods.

\begin{theorem}
 For the cyclic polling system with retrials and glue periods, the scaled steady-state joint queue length vector at an arbitrary time point in a visit period of station~$1$ satisfies
 \begin{equation}
(1-\rho)    \left(\begin{array}{cc}
  V_1^{(1q)}\\
        V_1^{(1o)} \\
        V_2^{(1)}\\
         \vdots \\
         V_N^{(1)}
        \end{array} \right)
\xrightarrow[d]{} \frac{b^{(2)}}{2b^{(1)}}\frac{1}{\delta} \left[  \left(\begin{array}{cc}
      (1-e^{-\nu_1 G_1})  \hat{u}^{(1)}_{1}  \\
     e^{-\nu_1 G_1}  \hat{u}^{(1)}_{1} \\
       \hat{u}^{(1)}_{2} \\
        \vdots \\
         \hat{u}^{(1)}_{N}
        \end{array} \right) + (1-e^{-\nu_1 G_1})  \hat{u}^{(1)}_{1} U \left(\begin{array}{cc}
       -1 \\
      \hat{\lambda}_1 E[B_1] \\
      \hat{\lambda}_2 E[B_1] \\
        \vdots \\
       \hat{\lambda}_N E[B_1]
        \end{array} \right)
        \right]
\Gamma(\alpha +1 , 1),~~~\text{when}~~ \rho \uparrow 1 .
\label{eq:randomvisit}
\end{equation}
\label{thm:randomvisit}
\end{theorem}
\proof{
We will use Equation (3.19) of \cite{Abidini16} to prove this. This equation states that the joint generating function, $R^{(i)}_{vi}(\underbar{$z$}_q,\underbar{$z$}_o)$, of the numbers of
customers in the queue and in the orbits at an arbitrary time point in a
visit period of $Q_i$ is given by
\begin{eqnarray*}
R^{(i)}_{vi}(\underbar{$z$}_q,\underbar{$z$}_o) &=& \frac{z_{iq} \left(\E [z_{iq}^{Y_i^{(iq)}}\left(\prod_{j=1,j\ne i}^{N}z_{jo}^{Y_j^{(i)}} \right)z_{io}^{Y_i^{(io)}}] - \E [\tilde{B_i}(\sum_{j=1}^{N} \lambda_j(1-z_{jo}))^{Y_i^{(iq)}}\left(\prod_{j=1,j\ne i}^{N}z_{jo}^{Y_j^{(i)}} \right)z_{io}^{Y_i^{(io)}}]\right)}{\E [Y_i^{(iq)}] \left(z_{iq} - \tilde{B_i}(\sum_{j=1}^{N} \lambda_j(1-z_{jo}))\right)} \nonumber \\
&& \times \frac{1-\tilde{B_i}(\sum_{j=1}^{N}\lambda_j(1-z_{jo}))} {\left(\sum_{j=1}^{N}\lambda_j(1-z_{jo})\right) \E [B_i]}.
\end{eqnarray*}
Evaluating the above generating function in the points $\underbar{$\tilde{z}$}_q = (z_{1q}^{(1-\rho)}, \ldots, z_{Nq}^{(1-\rho)})$ and $\underbar{$\tilde{z}$}_o = (z_{1o}^{(1-\rho)}, \ldots, z_{No}^{(1-\rho)})$, we get
 
\begin{eqnarray}
&&R^{(i)}_{vi}(\underbar{$\tilde{z}$}_q,\underbar{$\tilde{z}$}_o) = \nonumber \\
&& \frac{z_{iq}^{(1-\rho)} \left(\E [z_{iq}^{(1-\rho)Y_i^{(iq)}}\left(\prod_{j=1,j\ne i}^{N}z_{jo}^{(1-\rho)Y_j^{(i)}} \right)z_{io}^{(1-\rho)Y_i^{(io)}}] - \E [\tilde{B_i}(\sum_{j=1}^{N} \lambda_j(1-z_{jo}^{(1-\rho)}))^{Y_i^{(iq)}}\left(\prod_{j=1,j\ne i}^{N}z_{jo}^{(1-\rho)Y_j^{(i)}} \right)z_{io}^{(1-\rho)Y_i^{(io)}}]\right)}{\E [Y_i^{(iq)}] \left(z_{iq}^{(1-\rho)} - \tilde{B_i}(\sum_{j=1}^{N} \lambda_j(1-z_{jo}^{(1-\rho)}))\right)} \nonumber \\
&& \times \frac{1-\tilde{B_i}(\sum_{j=1}^{N}\lambda_j(1-z_{jo}^{(1-\rho)}))} {\left(\sum_{j=1}^{N}\lambda_j(1-z_{jo}^{(1-\rho)})\right) \E [B_i]}.
\label{copied_1}
\end{eqnarray}

This equation has two terms. The first term expresses the generating function of the number of customers in the system at the start of the service of the customer who is currently in service. The second term is the generating function of the number of customers that arrived during
the past service period of the customer who is currently in service.
As $\rho \uparrow 1$, this second term satisfies
\begin{eqnarray}
 \lim_{\rho \uparrow 1} \frac{1-\tilde{B_i}(\sum_{j=1}^{N}\lambda_j(1-z_{jo}^{(1-\rho)}))} {\left(\sum_{j=1}^{N}\lambda_j(1-z_{jo}^{(1-\rho)})\right) \E [B_i]} &=& \lim_{\rho \uparrow 1} \frac{1-\E[e^{-(\sum_{j=1}^{N}\lambda_j(1-z_{jo}^{(1-\rho)}))B_i}]} {\left(\sum_{j=1}^{N}\lambda_j(1-z_{jo}^{(1-\rho)})\right) \E [B_i]} \nonumber\\
 &=& \lim_{\rho \uparrow 1} \frac{\E[B_i e^{-(\sum_{j=1}^{N}\lambda_j(1-z_{jo}^{(1-\rho)}))B_i}]} { \E [B_i]} \nonumber\\
 &=& \frac{\E [B_i]}{\E [B_i]} = 1,
 \label{LSTRes}
\end{eqnarray}
where the second equality follows from l'H\^opital's rule.
Equation \eqref{LSTRes} expresses the fact that the scaled vector of number of customers arriving at the different stations during a past service time tends to $0$, and hence its generating function tends to $1$, as $\rho \uparrow 1$.

Before taking the limit $\rho \uparrow 1$ in \eqref{copied_1} we first look at
$ \lim_{\rho \uparrow 1} \left(\tilde{B_i}\left(\sum_{j=1}^{N} \lambda_j\left( 1-z_{jo}^{(1-\rho)}\right)\right)\right)^{\frac{1}{1-\rho}}$.
As we mentioned earlier when we scale $\rho \uparrow 1$ we scale the system such that only the arrival rates increase and the service times remain the same. So we can write $\lambda_j = \rho \hat{\lambda}_j$ where $\hat{\lambda}_j$ is fixed and independent
of $\rho$, for all $j= 1,\ldots, N.$ Therefore we have
\begin{eqnarray}
\lim_{\rho \uparrow 1} \left(\tilde{B_i}\left(\rho\sum_{j=1}^{N} \hat{\lambda}_j\left( 1-z_{jo}^{(1-\rho)}\right)\right)\right)^{\frac{1}{1-\rho}}
&=& e^{\lim_{\rho \uparrow 1} \frac{\ln \left(\tilde{B_i}\left(\rho\sum_{j=1}^{N} \hat{\lambda}_j     \left( 1-z_{jo}^{(1-\rho)}\right)\right)\right)}{1-\rho}} \nonumber\\
 &=& e^{\lim_{\rho \uparrow 1} \frac{\left(\rho \sum_{j=1}^{N} \hat{\lambda}_j   z_{jo}^{(1-\rho)} \ln z_{jo}  +  \sum_{j=1}^{N} \hat{\lambda}_j \left(1-z_{jo}^{(1-\rho)}\right)i\right)\tilde{B_i}^{\prime}\left(\rho\sum_{j=1}^{N}      \hat{\lambda}_j     \left( 1-z_{jo}^{(1-\rho)}\right)\right) }{-\tilde{B_i}\left(\rho\sum_{j=1}^{N}      \hat{\lambda}_j     \left( 1-z_{jo}^{(1-\rho)}\right)\right) }} \nonumber \\
  &=& e^{ \E[B_i] \sum_{j=1}^{N} \hat{\lambda}_j   \ln z_{jo} } =  e^{  \sum_{j=1}^{N}    \ln z_{jo}^{\E[B_i] \hat{\lambda}_j }} = \prod_{j=1}^{N}z_{jo}^{ \E[B_i] \hat{\lambda}_j} ,
\label{LSTHT}
\end{eqnarray}
where the second equality follows from l'H\^opital's rule.

Next we evaluate the following limit, which is related to the denominator of the first term in \eqref{copied_1}
\begin{eqnarray}
\lim_{\rho \uparrow 1} \frac{z_{iq}^{(1-\rho)} - \tilde{B_i}(\sum_{j=1}^{N} \lambda_j(1-z_{jo}^{(1-\rho)}))}{1 - \rho} &=& \lim_{\rho \uparrow 1} \frac{z_{iq}^{(1-\rho)}-1}{1 - \rho} +  \lim_{\rho \uparrow 1} \frac{1 - \E\left[e^{-\rho \sum_{j=1}^{N}  \hat{\lambda}_j \left(1-z_{jo}^{(1-\rho)}\right)B_i}\right]}{1 - \rho} \nonumber \\
&=&  \lim_{\rho \uparrow 1}z_{iq}^{(1-\rho)} \ln z_{iq}\nonumber\\
&-& \lim_{\rho \uparrow 1} \E\left[\left( \rho B_i \sum_{j=1}^{N} \hat{\lambda}_j   z_{jo}^{(1-\rho)} \ln z_{jo}  +  B_i\sum_{j=1}^{N} \hat{\lambda}_j \left(1-z_{jo}^{(1-\rho)}  \right)\right)e^{-\rho \sum_{j=1}^{N}  \hat{\lambda}_j \left(1-z_{jo}^{(1-\rho)}\right)B_i} \right] \nonumber \\
&=& \ln z_{iq}-  \sum_{j=1}^{N} \hat{\lambda}_j  \E [B_i] \ln z_{jo} = \ln \left( z_{iq}    \prod_{j=1}^{N}z_{jo}^{ - \E[B_i] \hat{\lambda}_j} \right),
\label{LSTdeno}
\end{eqnarray}
where the first equality uses the fact that $\lambda_j = \rho \hat{\lambda}_j$ and the second equality follows from l'H\^opital's rule.

We know that $\lim_{\rho \uparrow 1}  z_{iq}^{(1-\rho)} =1$. Substituting this along with \eqref{LSTRes},\eqref{LSTHT} and \eqref{LSTdeno} in \eqref{copied_1} while we take $\rho \uparrow 1$, we get
\begin{eqnarray}
 \lim_{\rho \uparrow 1}&&  R^{(i)}_{vi}(\underbar{$\tilde{z}$}_q,\underbar{$\tilde{z}$}_o) = \\
 &&\lim_{\rho \uparrow 1}  \frac{\E [z_{iq}^{(1-\rho)Y_i^{(iq)}}\left(\prod_{j=1,j\ne i}^{N}z_{jo}^{(1-\rho)Y_j^{(i)}} \right)z_{io}^{(1-\rho)Y_i^{(io)}}] - \E [\prod_{j=1}^{N}z_{jo}^{ \E[B_i] \hat{\lambda}_j (1-\rho)Y_i^{(iq)}}\left(\prod_{j=1,j\ne i}^{N}z_{jo}^{(1-\rho)Y_j^{(i)}} \right)z_{io}^{(1-\rho)Y_i^{(io)}}]}{\E [(1-\rho) Y_i^{(iq)}]  \ln \left( z_{iq}    \prod_{j=1}^{N}z_{jo}^{ - \E[B_i] \hat{\lambda}_j} \right)} .
\nonumber
\label{eq:metastate}
 \end{eqnarray}
Consider the following notation
\begin{eqnarray*}
\kappa &:=& \frac{b^{(2)}}{2b^{(1)}}\frac{1}{\delta}  (1-e^{-\nu_1 G_1})  \hat{u}^{(1)}_{1} =  \frac{b^{(2)}}{2b^{(1)}}\frac{\hat{\lambda}_1}{\delta}  \\
\kappa_1 &:=& \frac{b^{(2)}}{2b^{(1)}}\frac{1}{\delta}  e^{-\nu_1 G_1}  \hat{u}^{(1)}_{1} \\
\kappa_i &:=& \frac{b^{(2)}}{2b^{(1)}}\frac{1}{\delta} \hat{u}^{(1)}_{i}, ~~~~~~~~~~\forall i = 2, \ldots, N.
\end{eqnarray*}
Using the above notation in \eqref{eq:startvisit} and substituting it in \eqref{eq:metastate} we have,
\begin{eqnarray*}
 \lim_{\rho \uparrow 1}  R^{(1)}_{vi}(\underbar{$\tilde{z}$}_q,\underbar{$\tilde{z}$}_o) &=&   \frac{\E [z_{1q}^{\kappa \Gamma(\alpha , 1)}\prod_{j=1}^{N}z_{jo}^{\kappa_j \Gamma(\alpha , 1)} ] - \E [\prod_{j=1}^{N}z_{jo}^{ \E[B_1] \hat{\lambda}_j \kappa \Gamma(\alpha , 1)}\prod_{j=1}^{N}z_{jo}^{\kappa_j \Gamma(\alpha , 1)} ]}{\E [\kappa \Gamma(\alpha , 1)]  \ln \left( z_{1q}    \prod_{j=1}^{N}z_{jo}^{ - \E[B_1] \hat{\lambda}_j} \right)}  \nonumber \\
 &=& \frac{ \E [\left( z_{1q}^{\kappa} \prod_{j=1}^{N} z_{jo}^{\kappa_j} \right)^{\Gamma(\alpha , 1)} ] - \E [\left( \prod_{j=1}^{N}z_{jo}^{ \E[B_1] \hat{\lambda}_j \kappa + \kappa_j}\right)^{ \Gamma(\alpha , 1)} ] } {\kappa \alpha  \ln \left( z_{1q}    \prod_{j=1}^{N}z_{jo}^{ - \E[B_1] \hat{\lambda}_j} \right)}.
 \end{eqnarray*}
Now we introduce the following notation to change our generating function into an LST,
\begin{eqnarray*}
 s &:=& - \ln z_{1q} \\
 s_i &:=& - \ln z_{io}, ~~~~\forall i=1,\ldots,N.
\end{eqnarray*}
Then we have that the joint LST of the scaled steady-state joint queue length vector, of the queue of station~$1$ and the orbits at all the stations, during an arbitrary time in the visit period of station~$1$ is
\begin{eqnarray}
 \lim_{\rho \uparrow 1}  R^{(1)}_{vi}(\underbar{$\tilde{z}$}_q,\underbar{$\tilde{z}$}_o) &=& \frac{ \E [\left( e^{-s \kappa} \prod_{j=1}^{N} e^{-s_j \kappa_j} \right)^{\Gamma(\alpha , 1)} ] - \E [\left( \prod_{j=1}^{N}e^{-s_j \left(\E[B_1] \hat{\lambda}_j \kappa + \kappa_j\right)}\right)^{ \Gamma(\alpha , 1)} ] } { \alpha \kappa    \ln \left( e^{-s}    \prod_{j=1}^{N}e^{ s_j \E[B_1] \hat{\lambda}_j} \right)}  \nonumber \\
 &=& \frac{ \E [ e^{- \left(s \kappa + \sum_{j} s_j \kappa_j\right) \Gamma(\alpha , 1)} ] - \E [e^{ - \sum_{j=1}^{N}  s_j \left(\E[B_1] \hat{\lambda}_j \kappa + \kappa_j \right) \Gamma(\alpha , 1)} ]  } {\alpha \kappa  \left(-s+ \E[B_1] \sum_{j}  \hat{\lambda}_j s_j \right)}  \nonumber \\
 &=& \frac{\left( \frac{1}{1 +s \kappa + \sum_{j} s_j \kappa_j } \right)^{\alpha} - \left(\frac{1}{1 + \sum_{j=1}^{N}  s_j \left(\E[B_1] \hat{\lambda}_j \kappa + \kappa_j \right)} \right)^{\alpha}}{\alpha \kappa  \left(-s+ \E[B_1] \sum_{j}  \hat{\lambda}_j s_j \right)} \nonumber \\
 &=& \E [ e^{- \left(s \kappa + \sum_{j} s_j \kappa_j +    \left(-s \kappa+  \sum_{j}  \hat{\lambda}_j \E[B_1] \kappa s_j \right)  U \right) \Gamma(\alpha+1 , 1)} ],
 \label{eq:jointlst}
 \end{eqnarray}
where $U$ is a standard uniform random variable and the last equality follows from the expression
\[
 \E[e^{-(a+bU) \Gamma(\alpha+1,1)}] = \frac{\left(\frac{1}{1+a}\right)^{\alpha}-\left(\frac{1}{1+(a+b)}\right)^{\alpha}}{\alpha b}.
\]
Now we substitute $s = - \ln z_{1q}$ and $s_i = - \ln z_{io}$ back in \eqref{eq:jointlst} to get for the joint generating function
\begin{eqnarray}
 \lim_{\rho \uparrow 1}  R^{(1)}_{vi}(\underbar{$\tilde{z}$}_q,\underbar{$\tilde{z}$}_o) &=&  \E \left[\left(z_{1q}^{\kappa} \prod_{j=1}^{N} z_{jo}^{\kappa_j} \left(z_{1q}^{-1} \prod_{j=1}^{N} z_{jo}^{\hat{\lambda}_j \E[B_1]}  \right)^{ \kappa U} \right)^{\Gamma(\alpha+1 , 1)} \right].
 \label{eq:alpha1}
\end{eqnarray}
Let $V_{iq}^{(i)}$ and $V_{io}^{(i)}$ be the number of customers in the queue and orbit of station~$i$, and $V_{j}^{(i)}$ be the number of customers of type $j \ne i$, at an arbitrary point in time during a visit period of station~$i$, for all $i=1,\ldots,N$. Then from
\eqref{eq:alpha1} we have,
\begin{equation}
  (1-\rho)    \left(\begin{array}{cc}
  V_1^{(1q)}\\
        V_1^{(1o)} \\
        V_2^{(1)}\\
         \vdots \\
         V_N^{(1)}
        \end{array} \right)
\xrightarrow[d]{} \left[  \left(\begin{array}{cc}
       \kappa\\
     \kappa_1 \\
       \kappa_2 \\
        \vdots \\
         \kappa_{N}
        \end{array} \right) + \kappa   U \left(\begin{array}{cc}
       -1 \\
      \hat{\lambda}_1 E[B_1] \\
        \hat{\lambda}_2 E[B_1] \\
        \vdots \\
         \hat{\lambda}_N E[B_1]
        \end{array} \right)
       \right]
\Gamma(\alpha +1 , 1),~~~\text{when}~~ \rho \uparrow 1 .
\end{equation}
Therefore, the scaled steady-state joint queue length vector, in the queue of station~$1$ and the orbits of all stations, at an arbitrary time during a visit period of station~$1$, as $\rho \uparrow 1$, satisfies
\begin{equation*}
(1-\rho)    \left(\begin{array}{cc}
  V_1^{(1q)}\\
        V_1^{(1o)} \\
        V_2^{(1)}\\
         \vdots \\
         V_N^{(1)}
        \end{array} \right)
\xrightarrow[d]{} \frac{b^{(2)}}{2b^{(1)}}\frac{1}{\delta} \left[  \left(\begin{array}{cc}
      (1-e^{-\nu_1 G_1})  \hat{u}^{(1)}_{1}  \\
     e^{-\nu_1 G_1}  \hat{u}^{(1)}_{1} \\
       \hat{u}^{(1)}_{2} \\
        \vdots \\
         \hat{u}^{(1)}_{N}
        \end{array} \right) +  U \left(\begin{array}{cc}
       -\hat{\lambda}_1 \\
      \hat{\lambda}_1 \hat{\rho}_1 \\
      \hat{\lambda}_2 \hat{\rho}_1 \\
        \vdots \\
       \hat{\lambda}_N \hat{\rho}_1
        \end{array} \right)
        \right]
\Gamma(\alpha +1 , 1),~~~\text{when}~~ \rho \uparrow 1 .
\end{equation*} \qed}

An intuitive argument for Theorem \ref{thm:randomvisit} can be given in the following way.
Since, under heavy traffic, the scaled number of customers which are in the queue of station~$1$ at the start of an arbitrary visit period is
gamma distributed, $ \kappa \Gamma(\alpha , 1)$, also the scaled length
of an arbitrary visit period is gamma distributed, $ \kappa \E[B_1] \Gamma(\alpha , 1)$. Therefore if we choose an arbitrary time point in 
a visit period of station~$1$ the scaled length of that special visit period is distributed as  $ \kappa \E[B_1] \Gamma(\alpha +1 , 1)$,
where $\kappa= \frac{b^{(2)}}{2 b^{(1)}}\frac{1}{\delta }(1-e^{-\nu_1 G_1}) \hat{u}^{(1)}_{1}$.  Since this is a special interval
which we are looking at, the scaled steady-state joint queue length vector at the start of this visit period satisfies

 \begin{equation}
  (1-\rho)    \left(\begin{array}{cc}
  \breve{Y}_1^{(1q)}\\
        \breve{Y}_1^{(1o)} \\
        \breve{Y}_2^{(1)}\\
         \vdots \\
         \breve{Y}_N^{(1)}
        \end{array} \right)
\xrightarrow[d]{}\frac{b^{(2)}}{2b^{(1)}}\frac{1}{\delta}   \left(\begin{array}{cc}
       (1-e^{-\nu_1 G_1})  \hat{u}^{(1)}_{1} \\
     e^{-\nu_1 G_1}  \hat{u}^{(1)}_{1} \\
       \hat{u}^{(1)}_{2} \\
        \vdots \\
         \hat{u}^{(1)}_{N}
        \end{array} \right)
\Gamma(\alpha +1 , 1),~~~\text{when}~~ \rho \uparrow 1 .
\label{eq:specialstart}
 \end{equation}
 At the arbitrary point in time, we have $ \kappa U \Gamma(\alpha +1 , 1)$ customers served, which means there are $ J_j^{(1)} = \sum_{k=1}^{ \kappa U \Gamma(\alpha +1 , 1)} (L_j^{(1)})_k$ new arrivals
 of type~$j$ customers during that period. Note that as $\rho \uparrow 1$, $ J_j^{(1)} \rightarrow \infty$ , therefore the new arrivals during the past service time of the customer in service can be neglected.
 Using the same arguments as in Theorem \ref{startofswitch} we can say that the limiting scaled distribution of the new number of customers
 of type~$j$ at an arbitrary point in time during the visit of station~$i$ can be given as
 \[
\frac{ J_j^{(1)}}{U  \breve{Y}_1^{(1q)}} \rightarrow \lambda_j \E[B_1] .
 \]
Therefore the scaled steady-state joint queue length vector at an arbitrary point in time during the visit of station~$1$ as $\rho \uparrow 1$ satisfies
\[(1-\rho)    \left(\begin{array}{cc}
  V_1^{(1q)}\\
        V_1^{(1o)} \\
        V_2^{(1)}\\
         \vdots \\
         V_N^{(1)}
        \end{array} \right)
\xrightarrow[d]{}
 \frac{b^{(2)}}{2b^{(1)}}\frac{1}{\delta} \left[  \left(\begin{array}{cc}
       0 \\
     e^{-\nu_1 G_1}  \hat{u}^{(1)}_{1} \\
       \hat{u}^{(1)}_{2} \\
        \vdots \\
         \hat{u}^{(1)}_{N}
        \end{array} \right) + (1-e^{-\nu_1 G_1})  \hat{u}^{(1)}_{1}  \left(\begin{array}{cc}
      1- U \\
     U \hat{\lambda}_1 E[B_1] \\
        U \hat{\lambda}_2 E[B_1] \\
        \vdots \\
         U \hat{\lambda}_N E[B_1]
        \end{array} \right)
        \right]
\Gamma(\alpha +1 , 1),~~~\text{when}~~ \rho \uparrow 1,
\]
which is equivalent to \eqref{eq:randomvisit}.

In \eqref{eq:randomvisit} we have given the scaled steady-state joint queue length vector of customers of each type at an arbitrary point in time during the visit period of station~$1$ when $\rho \uparrow 1$. Using Remark \ref{rem:vector}
we can extend this to an arbitrary point in time during the visit period of a given station~$i$. This can be written as
\begin{equation}
(1-\rho)    \left(\begin{array}{cc}
        V_1^{(i)}\\
         \vdots \\
         V_{i-1}^{(i)}\\
V_i^{(iq)}\\
        V_i^{(io)} \\
                V_{i+1}^{(i)}\\
         \vdots \\
         V_N^{(i)}
        \end{array} \right)
\xrightarrow[d]{} \frac{b^{(2)}}{2b^{(1)}}\frac{1}{\delta} \left[  \left(\begin{array}{cc}
        \hat{u}^{(i)}_{1} \\
        \vdots \\
         \hat{u}^{(i)}_{i-1}\\
      (1-e^{-\nu_i G_i})  \hat{u}^{(i)}_{i}  \\
     e^{-\nu_i G_i}  \hat{u}^{(i)}_{i} \\
       \hat{u}^{(i)}_{i+1} \\
        \vdots \\
         \hat{u}^{(i)}_{N}
        \end{array} \right) + U \left(\begin{array}{cc}
       \hat{\lambda}_1 \hat{\rho}_i \\
        \vdots \\
       \hat{\lambda}_{i-1} \hat{\rho}_i\\
       -\hat{\lambda}_i \\
      \hat{\lambda}_i \hat{\rho}_i \\
      \hat{\lambda}_{i+1} \hat{\rho}_i \\
        \vdots \\
       \hat{\lambda}_N \hat{\rho}_i
        \end{array} \right)
        \right]
\Gamma(\alpha +1 , 1),~~~\text{when}~~ \rho \uparrow 1 .
\label{eq:randomvisit_i}
\end{equation}
Due to the observation that in heavy traffic, the server resides in a visit period for 100\% of the time, \eqref{eq:randomvisit_i} leads to the following theorem.
\begin{theorem}\label{thm:final}
In a cyclic polling system with retrials and glue periods, the scaled steady-state joint queue length vector at an arbitrary time point, with $L^{(iq)}$ and $L^{(io)}$ representing the number in queue and in orbit at station $i$ respectively for all $i=1,\ldots,N$,
satisfies
\begin{equation}
(1-\rho)    \left(\begin{array}{cc}
        L^{(1q)}\\
        \vdots \\
         L^{(Nq)}\\
         L^{(1o)}\\
	 \vdots\\
	 L^{(No)}
        \end{array} \right)
\xrightarrow[d]{}\frac{b^{(2)}}{2b^{(1)}}\frac{1}{\delta} \underbar{P}~ \Gamma(\alpha +1 , 1) ~~~\text{when}~~ \rho \uparrow 1,
\label{eq:thmrandompoint}
\end{equation}
 where $\underbar{P} = \underbar{P}_i$ with probability $\hat{\rho}_i$ and
\[
 \underbar{P}_i = \left(\begin{array}{cc}
 0\\

        \vdots \\
        0\\

      (1-e^{-\nu_i G_i})  \hat{u}^{(i)}_{i}  \\
     0\\
     \vdots \\
     0\\
       \hat{u}^{(i)}_{1} \\
       \vdots\\
        \hat{u}^{(i)}_{i-1}\\
         e^{-\nu_i G_i}  \hat{u}^{(i)}_{i} \\
       \hat{u}^{(i)}_{i+1} \\
        \vdots \\
         \hat{u}^{(i)}_{N}
        \end{array} \right) +  U \left(\begin{array}{cc}
        0\\
        \vdots \\
        0\\
        -\hat{\lambda}_i\\
        0\\
        \vdots\\
        0\\
       \hat{\lambda}_1 \hat{\rho}_i \\
        \vdots \\
       \hat{\lambda}_{i-1} \hat{\rho}_i\\
      \hat{\lambda}_i \hat{\rho}_i \\
      \hat{\lambda}_{i+1} \hat{\rho}_i \\
        \vdots \\
       \hat{\lambda}_N \hat{\rho}_i
        \end{array} \right)   .
\]
\end{theorem}
\proof{
As mentioned at the beginning of this section, when $\hat{\rho} \uparrow 1$ the system is in the visit periods with probability 1. Therefore
the limiting scaled joint queue length distribution at an arbitrary point in time can be given as
the limiting scaled joint queue length distribution in the visit period of station~$i$ with probability $\hat{\rho}_i.$ 
Now consider that the number of customers of type~$j$, at an arbitrary point in time during the visit period of station~$i$, in queue and orbit respectively is given by $V_j^{(iq)}$ and $V_j^{(io)}$. Then using
\eqref{eq:randomvisit_i} we can write
\begin{equation}
(1-\rho)    \left(\begin{array}{cc}
       V_1^{(iq)}\\
        \vdots \\
       V_{i-1}^{(iq)}\\
        V_i^{(iq)}\\
               V_{i+1}^{(iq)}\\
                \vdots \\
       V_N^{(iq)}\\
        V_1^{(io)}\\
        \vdots\\
         V_{i-1}^{(io)}\\
         V_{i}^{(io)}\\
         V_{i+1}^{(io)}\\
      \vdots \\
      V_N^{(io)}
        \end{array} \right)
\xrightarrow[d]{} \frac{b^{(2)}}{2b^{(1)}}\frac{1}{\delta} \left(\left(\begin{array}{cc}
 0\\

        \vdots \\
        0\\

      (1-e^{-\nu_i G_i})  \hat{u}^{(i)}_{i}  \\
     0\\
     \vdots \\
     0\\
       \hat{u}^{(i)}_{1} \\
       \vdots\\
        \hat{u}^{(i)}_{i-1}\\
         e^{-\nu_i G_i}  \hat{u}^{(i)}_{i} \\
       \hat{u}^{(i)}_{i+1} \\
        \vdots \\
         \hat{u}^{(i)}_{N}
        \end{array} \right) + U \left(\begin{array}{cc}
        0\\
        \vdots \\
        0\\
        -\hat{\lambda}_i\\
        0\\
        \vdots\\
        0\\
       \hat{\lambda}_1 \hat{\rho}_i \\
        \vdots \\
       \hat{\lambda}_{i-1} \hat{\rho}_i\\
      \hat{\lambda}_i \hat{\rho}_i \\
      \hat{\lambda}_{i+1} \hat{\rho}_i \\
        \vdots \\
       \hat{\lambda}_N \hat{\rho}_i
        \end{array} \right)\right)
\Gamma(\alpha +1 , 1),~~~\text{when}~~ \rho \uparrow 1 .
\end{equation}

This holds because $V_j^{(io)} = V_j^{(i)}$ and $V_j^{(iq)} = 0$ when $i \ne j.$
Therefore we know that the limiting scaled joint queue length distribution at an arbitrary point in time, with probability $\hat{\rho}_i$, can be given as
$  \frac{b^{(2)}}{2b^{(1)}}\frac{1}{\delta}  \underbar{P}_i\Gamma(\alpha +1 , 1)$. Hence we have
\begin{equation}
(1-\rho)    \left(\begin{array}{cc}
        L^{(1q)}\\
        \vdots \\
         L^{(Nq)}\\
           L^{(1o)}\\
           \vdots\\
        L^{(No)}
        \end{array} \right)
\xrightarrow[d]{}\frac{b^{(2)}}{2b^{(1)}}\frac{1}{\delta} \underbar{P}~ \Gamma(\alpha +1 , 1), ~~~\text{when}~~ \rho \uparrow 1,
\end{equation}
 where $\underbar{P} = \underbar{P}_i$ with probability $\hat{\rho}_i$.
\qed}


\begin{remark}

Note that under heavy-traffic the total scaled workload in the system satisfies the so-called
\emph{heavy-traffic averaging principle}. This principle, first found in \cite{CoffmanEtAl1,CoffmanEtAl2}
for a specific class of polling models, implies that the workload in each queue is emptied and refilled
at a rate that is much faster than the rate at which the total workload is changing. As a consequence the total workload
can be considered constant during the course of a cycle (represented by the gamma distribution), while the workloads
in the individual queues fluctuate like a fluid model. It is because of this that the queue length vector
in Theorem \ref{thm:final} also features a state-space collapse (cf.\ \cite{Reiman}): the limiting distribution of
the $2N$-dimensional scaled queue length vector is governed by just three distributions: the discrete distribution governing
\underbar{P}, the uniform distribution and the gamma distribution.

Therefore using Theorems \ref{thm:randomvisit} and \ref{thm:final}, and the fact that $\delta =\sum_i \E[B_i] \hat{u}^{(j)}_{i}$ 
for all $j = 1, \ldots, N,$ the scaled workload in the system at an arbitrary point in time is given by
\begin{equation*}
(1-\rho)\sum_{i=1}^N \E[B_i] (L^{(iq)}+L^{(io)})\xrightarrow[d]{} \frac{b^{(2)}}{2b^{(1)}}\Gamma(\alpha+1, 1),~~~\text{when}~~ \rho \uparrow 1 .
\end{equation*}

Since the above equation is independent of everything but the gamma distribution, the workload is the same at any arbitrary point in time during the cycle.
Hence the system agrees with the \emph{heavy-traffic averaging principle}. Note that in Theorems \ref{thm:startglue}, \ref{startofvisit} and \ref{startofswitch}
we have found the limiting distribution of the scaled number of customers at embedded time points. Extending the heavy traffic principle along with these theorems, we can say that the
scaled workload in an arbitrarily chosen cycle can be given as
\begin{equation*}
(1-\rho)\sum_{i	=1}^N \E[B_i] X_i^{(1)}\xrightarrow[d]{} \frac{b^{(2)}}{2b^{(1)}}\Gamma(\alpha, 1),~~~\text{when}~~ \rho \uparrow 1 .
\end{equation*}

We observe that the scaled workload in the two cases, arbitrary point in time and arbitrary cycle, have $\Gamma(\alpha+1, 1)$ and $\Gamma(\alpha, 1)$ distributions respectively.
This is because of a bias introduced in selection of an arbitrary point in time, i.e., an arbitrarily chosen point in time has a higher probability to be in a longer cycle than being in a shorter cycle. This
bias does not exist when we arbitrarily choose a cycle.

\end{remark}

\section{Approximations}

In the previous section we derived the heavy traffic limit of the scaled steady-state joint queue length vector
at an arbitrary point in time for the cyclic polling system with retrials and glue periods. We will use this result to derive
an approximation for the mean number of customers in the system for arbitrary system loads. This is done by using an interpolation between the heavy traffic result and a light traffic result, similar to what is described in Boon et al. \cite{Boon2011}.

\subsection{Approximate mean number of customers}
Consider the following approximation for the mean number of customers of type~$i$,
\begin{equation}\label{eq:approx}
 E[L_i] \approx \frac{c_0 + \rho c_1}{1-\rho}.
\end{equation}
The coefficients $c_0$ and $c_1$ are chosen in agreement with the light traffic and heavy traffic behaviour of $E[L_i]$. Clearly, when $\rho \downarrow 0$, also $E[L_i] \downarrow 0$.
Hence we choose $c_0=0$.
On the other hand, we have $ \lim_{\rho \uparrow 1} E[(1-\rho) L_i] = c_1$. Using \eqref{eq:thmrandompoint}, we get
\begin{eqnarray*}
 \lim_{\rho \uparrow 1} E[(1-\rho) L_i] &=&  \lim_{\rho \uparrow 1}\E[(1-\rho)( L^{(iq)}+ L^{(io)}) ] \\
&=& \E\left[ \frac{b^{(2)}}{2b^{(1)}}\frac{1}{\delta}  \left(\sum_{j=1}^{N} \hat{\rho}_j \left(\hat{u}_i^{(j)} +U \hat{\lambda}_i\hat{\rho}_j\right) - U \hat{\lambda}_i\hat{\rho}_i \right) \Gamma(\alpha +1 , 1) \right]\\
&=& \frac{b^{(2)}}{2b^{(1)}}\frac{1}{\delta} \left( \sum_{j=1}^{N} \hat{\rho}_j \left(\hat{u}_i^{(j)} + \frac{ \hat{\lambda}_i\hat{\rho}_j}{2}\right) - \frac{ \hat{\lambda}_i\hat{\rho}_i}{2} \right) (\alpha +1).
\end{eqnarray*}
Hence we choose $ c_1 = \frac{b^{(2)}}{2b^{(1)}}\frac{(\alpha +1)}{\delta} \left( \sum_{j=1}^{N} \hat{\rho}_j \left(\hat{u}_i^{(j)} + \frac{ \hat{\lambda}_i\hat{\rho}_j}{2}\right) - \frac{ \hat{\lambda}_i\hat{\rho}_i}{2} \right) $, and so the final approximation for $E[L_i]$ becomes, for arbitrary
$\rho \in (0,1)$,
\begin{eqnarray}
  E[L_i] \approx \frac{\rho}{(1-\rho)} \frac{b^{(2)}}{2b^{(1)}}\frac{(\alpha +1)}{\delta} \left( \sum_{j=1}^{N} \hat{\rho}_j \left(\hat{u}_i^{(j)} + \frac{ \hat{\lambda}_i\hat{\rho}_j}{2}\right) - \frac{ \hat{\lambda}_i\hat{\rho}_i}{2} \right) .
\label{eq:approxnumberi}
  \end{eqnarray}

\subsection{Numerical results}

In this section we will compare the above approximation with exact results. The exact results are obtained using the approach in \cite{Abidini16}.

Consider a five station polling system in which the service times are exponentially distributed with mean $\E[B_i] =1$ for all $i = 1, \ldots, 5$. The arrival processes are Poisson processes with rates $\lambda_1 = \rho\frac{1}{10},
\lambda_2 = \rho\frac{2}{10},~\lambda_3 = \rho\frac{3}{10},~\lambda_4 = \rho\frac{1}{10},~\lambda_5 = \rho\frac{3}{10}.$  The switch-over times from station~$i$ are exponentially distributed with mean
$\E[S_i] = 2,~3,~1,~5,~2$ for stations $i=1, \ldots, 5$. The durations of the deterministic glue periods are $G_i = 3,~1,~2,~1,~2$, and the exponential retrial rates are $\nu_i = 5,~1,~3,~2,~1$, for stations $i=1$ to $5$ respectively.
We plot the following for $\rho \in (0,1)$ and compare the approximation given in \eqref{eq:approxnumberi} with the values obtained using exact analysis.

\begin{multicols}{2}
      {\centering
      \includegraphics[width=0.5\textwidth]{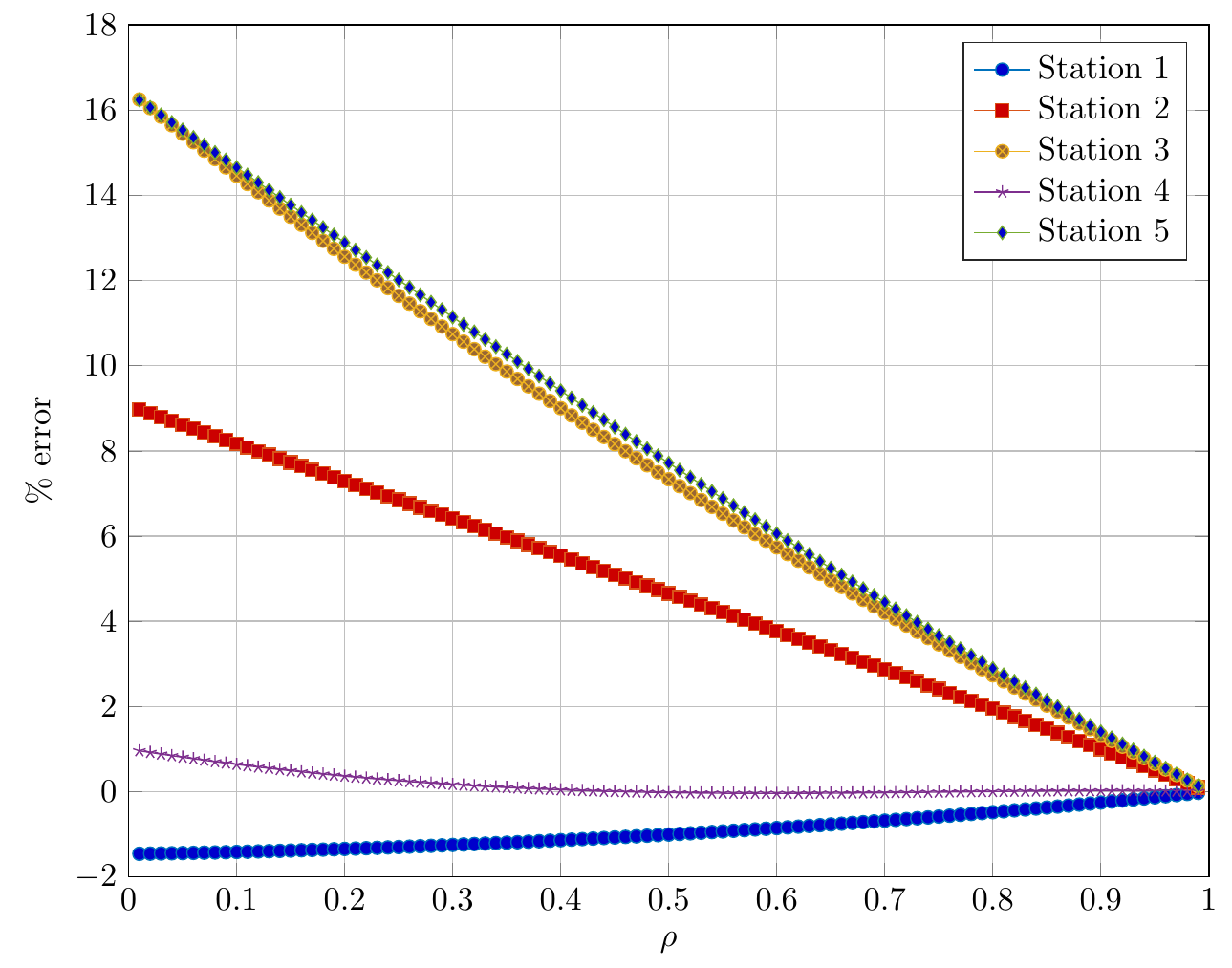}
        \captionof{figure}{$\%$ error for the number of customers in each station}
        \label{fig:gull1}}
{\centering
        \includegraphics[width=0.5\textwidth]{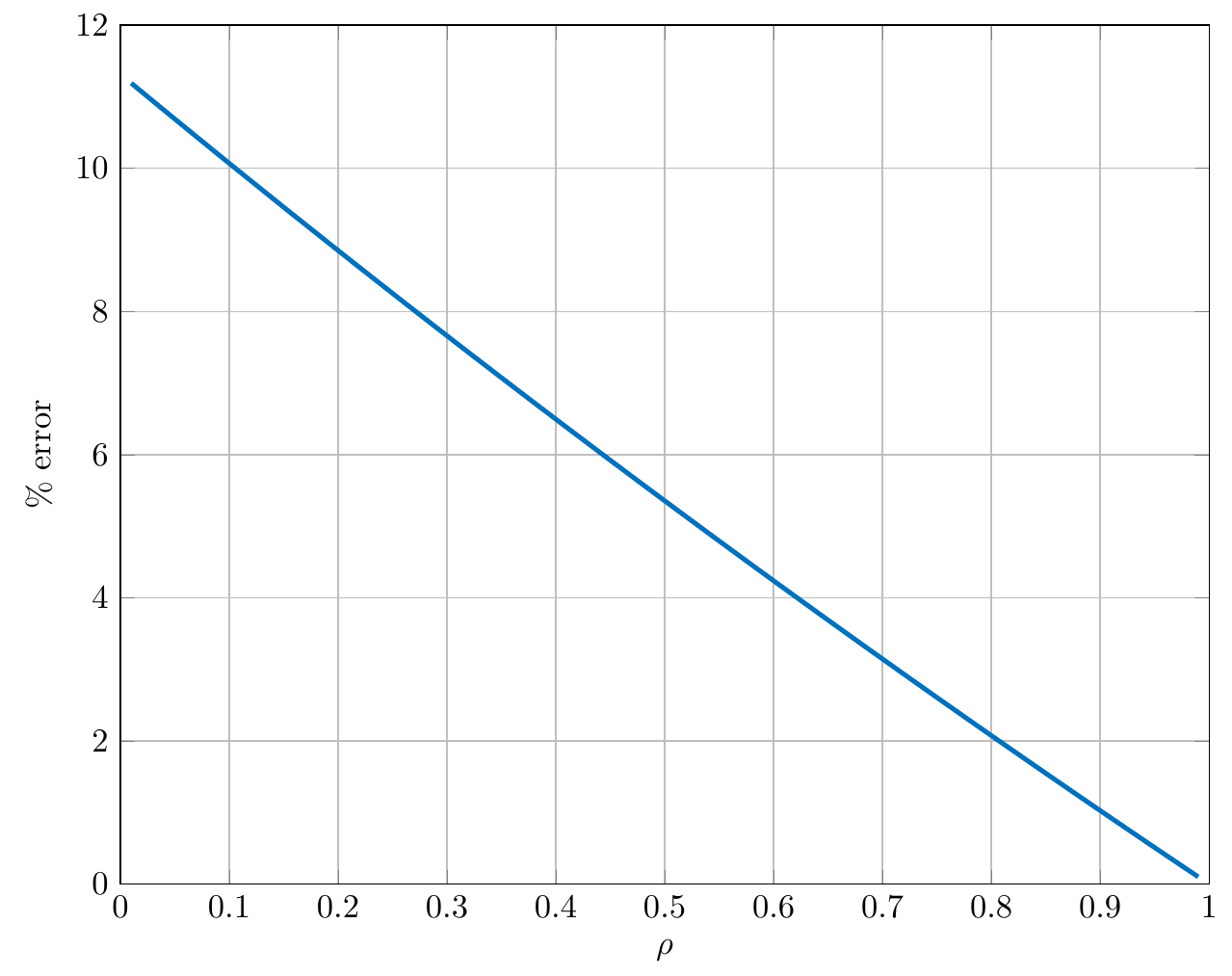}
        \captionof{figure}{$\%$ error for total number of customers}
      \label{fig:tiger1}}
\end{multicols}

In Figures \ref{fig:gull1} and \ref{fig:tiger1} we respectively plot the percentage error calculated as $ \%~\text{error}= \frac{\text{Approximate Value}-\text{Exact Value} }{\text{Exact Value}} \times 100$, for the mean
number of customers of each type and the total mean number of customers in the system. The error percentage is similar to that predicted in \cite{Boon2011}. The error is non-negligible for lower values of $\rho$, but it decreases quickly as $\rho$ increases. Consequently, for larger values of $\rho$, the approximation is accurate.

Based on this, we conclude that the heavy-traffic results as derived in this paper are very useful for deriving closed-form approximations for the queue length, especially as the systems under study (e.g.\ optical systems) typically run under a heavy workload (i.e., a large value of $\rho$). Nevertheless, to obtain better performance for small values of $\rho$, the current approximation as presented here can be refined by e.g.\ computing theoretical values of $\frac{d}{d\rho}\E{L_i}|_{\rho = 0}$ and incorporating that information in \eqref{eq:approx} as explained in \cite{Boon2011}. Furthermore, approximations for the mean queue length as mentioned here can be extended to approximations for the complete queue length distributions of the polling systems with glue periods and retrials in the spirit of \cite{DorsmanEtAl}. These extensions, however, are beyond the scope of this paper.


\subsection*{Acknowledgement}
\noindent
The authors wish to thank Marko Boon and Onno Boxma for fruitful discussions.
The research is supported by the IAP program BESTCOM funded by the Belgian government, and by the Gravity program NETWORKS funded
by the Dutch government. Part of the research of the second author was performed while he was affiliated with Leiden University.

\end{document}